\documentclass[12pt]{article}

\usepackage{amsmath, amssymb, amsthm, epsfig, amsfonts}

\def\al{\alpha}
\def\b{\beta}

\def\vp{\varphi}
\def\S{\Sigma}

\def\La{\Lambda}
\def\ck{\check}
\def\wh{\widehat}
\def\wt{\widetilde}

\newfont{\msbm}{msbm10}
\newfont{\msbms}{msbm6}

\def\N{{\mathbb N}}
\def\P{{\mathbb P}}
\def\C{{\mathbb C}}
\def\R{{\mathbb R}}
\def\Z{{\mathbb Z}}
\def\HH{{\mathbb H}}
\def\h{\mathfrak{h}}
\def\H{{\cal H}}
\def\M{{\cal M}}
\def\I{{\cal I}}

\def\bi{\begin{itemize}}
\def\ei{\end{itemize}}
\def\la{\label}
\def\nn{\nonumber}

\def\l2{L^2(T^g, {d^g \vp})}

\def\j0{J^0(\S)}

\def\a0{{\cal A}_0^{SU(r)}/{\cal G}}

\def\u1r{U(1)^{r-1}}

\newtheorem{theorem}{Theorem}
        \newtheorem{proposition}{Proposition}
        \newtheorem{lemma}{Lemma}
        \newtheorem{corollary}{Corollary}
        \newtheorem{definition}{Definition}
        \newtheorem{problem}{Problem}
        \newtheorem{remark}{Remark} 

\newcommand{\ba}{\begin{eqnarray}}
\newcommand{\ea}{\end{eqnarray}}

\title{Coherent State Transforms and Vector Bundles on Elliptic Curves}

\author{Carlos A. Florentino, Jos\'e M. Mour\~ao and Jo\~ao P. Nunes\\
\small{cfloren, jmourao, jpnunes@math.ist.utl.pt}\\
\small{Department of Mathematics, Instituto Superior T\'ecnico},\\ 
\small{Av. Rovisco Pais, 1049-001 Lisboa, Portugal}}

\begin{document}

\maketitle

\vskip 0.5cm

\vspace{.5cm}

\abstract{
We extend the coherent state transform (CST) of Hall to the context 
of the moduli spaces of semistable holomorphic  vector bundles with fixed determinant 
over elliptic curves. We show that by applying the CST to appropriate distributions, 
we obtain the space of level $k$, rank $n$ and genus one non-abelian 
theta functions with the unitarity of the CST transform being preserved. Furthermore, the  
shift $k \rightarrow k+n$ appears in a natural 
way in this finite-dimensional framework.
}

\newpage

\tableofcontents

\newpage


\section{Introduction}

In \cite{Ha1} Hall proposed a generalization of 
the Segal-Bargmann or coherent state transform (CST) \cite{Se1,Se2,Ba}   in which $\R^n$ is replaced by an arbitrary
compact connected Lie group $K$ and $\C^n$ by
the complexification $K_\C$ of $K$.
This Segal-Bargmann-Hall CST was further generalized
to gauge theories with applications to gravity in the context of Ashtekar 
variables
in \cite{ALMMT} and to Yang-Mills theories in two 
space-time dimensions in \cite{DH}. For reviews and further developments see
\cite{Ha2,Ha3} and \cite{Th}.

In the present paper we continue the project started in \cite{FMN}
of the application of CST techniques to the study of theta functions.

For a Riemann surface $X$ of genus $g$, 
$K_\C$-theta functions are sections of holomorphic line
bundles ${\cal L}$ over the moduli space of semistable $K_\C$-bundles on
$X$
\ba
{\cal L} & \to & {\cal M}_{K_\C}(X) \nonumber \\
\theta &   \in & H^0({\cal M}_{K_\C}(X), {\cal L}) \nonumber .
\ea
The study of theta functions motivates considerable interest both
from the mathematical and physical points of view. 
In physics, the spaces  
 $H^0({\cal M}_{K_\C}(X), {\cal L})$ correspond both to spaces of
conformal blocks in WZW conformal field theories and to Hilbert spaces
of states of Chern-Simons theories. It was in this context that the Verlinde 
formula for the dimensions of these vector spaces of holomorphic sections 
was discovered, see for example \cite{So} for a review.

In the case when $K_\C$ is $SL(n,\C)$,
this moduli space can be interpreted, via the well-known theorem of 
Narasimhan and Seshadri \cite{NS}, as the moduli space ${\cal M}_n(X):={\cal M}_{SL(n,\C)}(X)$ of semistable rank $n$
vector bundles with trivial determinant over $X$, and the corresponding 
non-abelian theta functions were already subject of study by Weil. 
The conformal blocks are then represented by holomorphic sections of powers of  
determinant line bundles over ${\cal M}_n(X)$ which is also the moduli space of flat 
$SU(n)$-connections on $X$. 
These non-abelian theta functions have been widely studied since the
nineties, mainly from the point of view of algebraic geometry. 
However, an analytic theory of these functions is not yet fully developed, 
and there are many open questions related to them \cite{Bea,Fa}. 

In this work we will consider the case when $X$ is an elliptic curve, $X = X_\tau, {\rm Im}\tau >0$. 
Non-abelian theta functions of level $k$ on elliptic curves have been studied 
mainly from two points of view. Expressions for orthonormal frames for these theta functions 
have been obtained by physical (formal)
functional integral methods in conformal field theory, e.g. in  \cite{Ber,EMSS,FG,G}.
These expressions were also obtained from infinite dimensional geometric quantization and 
symplectic reduction of affine spaces in the context of Chern-Simons theory \cite{AdPW}. In both approaches, 
one observes a shift in the level $k \rightarrow k+h$, with $h$ being the dual Coxeter number of $K$ ($h=n$ for $SU(n)$), 
which 
arises from a regularization of infinite determinants of differential operators on bundles over $X$.   
We will show that by extending the CST of Hall for $SU(n)$ to appropriate finite-dimensional spaces 
of distributions, we obtain the spaces of non-abelian theta functions and that the averaged heat kernel 
measure descends to a hermitean structure on ${\cal L}^k$ making the CST transform unitary, with the 
correct shift of level $k \rightarrow k+n$.  

More precisely, to obtain the relation between the CST for $SU(n)$ and
genus one non-abelian theta functions 
we adopt a strategy similar to that of section 4 of \cite{FMN}:
\begin{itemize}

\item[1)]
In propositions \ref{p41} and \ref{p42} of subsection \ref{vemdo4.2} we consider  
the extension of the CST
to complexified Laplacians and then to the space of distributions
$C^\infty (SU(n))^\prime$. 
By taking the second step we of course lose unitarity of the CST. The restriction of the CST 
for a simple compact group $K$ to $Ad_K$-invariant functions and distributions is intimately related 
with the CST for the maximal torus $T\subset K$, but for Weyl anti-invariant functions and distributions as 
we show in theorems \ref{tb1} and \ref{pcar}. This simple fact will play a crucial role and it is behind the 
success of the CST in reproducing the shift $k \rightarrow k+h$.

\item[2)] The relation between holomorphic functions on $SL(n,\C)$ (obtained from this extended CST) and 
sections of line bundles over ${\cal M}_n(X_\tau)$,  is provided by pull-back via the Schottky map, 
see section \ref{s3} (specially equation (\ref{nuncamaisacaba}) and proposition \ref{semlabel}), 
$$
S: SL(n,\C) \rightarrow {\cal M}_n(X_\tau).
$$
In analogy with the abelian case, this gives a description of non-abelian theta functions in genus 
one as holomorphic functions on a space whose complex structure is canonical (independent of $\tau$). 

\item[3)]
The theorems of subsection \ref{s42} show that the restriction of
the CST to appropriate finite dimensional subspaces ${\cal F}_k$
of Ad-invariant distributions (see (\ref{510}) and (\ref{tasol})) leads to a vector bundle
over the Teichm\"uller space of genus one curves
\begin{equation}
\la{4221}
\wt \H_k \rightarrow {\cal T}_1= \{ \tau \in \C: {\rm Im}\tau >0\}
\end{equation}
which is
isomorphic to the 
vector bundle $\H_k  \rightarrow {\cal T}_1 $ of conformal blocks, see (\ref{4121}),
with simple unitary (see the point 4) below) isomorphism $\Phi_k$ given by
\begin{eqnarray}
\la{4220}
\Phi_{k,\tau}: \wt \H_{k,{\tau}} & \rightarrow & \H_{k,\tau}^+ \\
\nonumber
\Psi & \mapsto & \theta^+ = e^{\frac{||\rho||^2}{k+n} 
\pi i \tau}\frac{\sigma}{\theta_{\rho,n}^-}\Psi 
\end{eqnarray} 
where $\wt \H_{k,{\tau}} = \wt \H_{k|_{\tau}}$, 
$\H_{k,\tau}^+=\H_{k|_{\tau}}$ 
and $\Phi_{k,\tau} = \Phi_{k|_{\widetilde \H_{k,{\tau}}}}$. 
Here, $\theta^{\pm}$ 
denote Weyl
invariant or anti-invariant theta functions, $\rho$ is the Weyl vector
and $\sigma$ is the denominator in the Weyl character formula.
The image of ${\cal F}_k$ under the CST selects a
trivialization of the bundle of conformal blocks corresponding
to the frame (\ref{418}). This shows
that the CST on $SU(n)$ leads to the shift of level
$k \rightarrow k+n$ also obtained in conformal field
theory but with the help of 
infinite dimensional Feynman path integral methods.

\item[4)]
The Hall averaged heat kernel measure on $SL(n,\C)$
defines a hermitian structure (\ref{831}) on $\wt \H_k$
for which the unitarity of 
the extended CST is recovered.
\end{itemize}

This strategy provides an intrinsically finite-dimensional framework for 
the study of non-abelian theta functions in genus one.

The organization of the paper is the following. 
In section \ref{s2}, we extend the CST to class functions on compact 
Lie groups and extend the results to $Ad$-invariant distributions on
the group.
In section \ref{s22}, we extend the results of 
\cite{FMN} to abelian varieties with a general polarization. 
In the following sections, we restrict ourselves 
to the context of elliptic curves and $K=SU(n)$
($K_\C=SL(n,\C)$) corresponding to the moduli space ${\cal M}_n(X_\tau)$.
(The corresponding results for other Lie groups should also apply.)
In section \ref{s3}, we describe the Schottky map from the space of Schottky representations 
of $\pi_1(X_\tau)$ in $K_\C$ to the moduli space $\M_n(X_\tau)$.   
By considering the Schottky map, we show in section \ref{s4} that the results in \cite{FMN} extend to non-abelian theta functions 
in genus 1.
Namely, applying the CST to appropriate finite dimensional
spaces of distributions ${\cal F}_k\subset C^\infty(SU(n))'$
yields the spaces of level $k$ non-abelian theta functions
(Theorem \ref{levico8}).
In this case however, as we show in Theorems \ref{levico8} and \ref{levico9}, 
not only the Hall averaged heat kernel measure
makes the CST unitary, but also the hermitean structure that it defines
on the corresponding vector bundle over the Teichm\"uller space is
the ``correct'' one in the sense that it contains the level shift
$k\mapsto k+n$. Surprisingly, the case of $SU(2)$ is special and is treated 
in subsection \ref{5.3}, where the main results are formulated in theorem \ref{maio7}.

Extensions of this work and applications to the moduli space of semistable 
vector bundles on higher genus curves will appear in \cite{FMNT}.


\section{Extensions of the Coherent State Transform}
\label{s2}
 

\subsection{Coherent State Transform for Lie Groups}
\label{s21a}

Let $K$ be a compact connected Lie group of rank $l$, 
$K_{\C}$ its complexification  (see \cite{Ho}) and let
$\rho_{t}, t>0,$ be the heat kernel 
for the Laplacian $\Delta$ on $K$ associated to an Ad-invariant 
inner product on its Lie algebra $Lie(K)$. 
If $\{X_i, \, i=1,...,\dim K\}$ is a corresponding  orthonormal 
basis for $Lie(K)$ viewed as the space 
of left-invariant vector fields on $K$, then $\Delta = 
\sum_{i=1}^{\dim K} X_iX_i$. 
As proved in \cite{Ha1}, $\rho_{t}$ has a unique analytic 
continuation to $K_{\C}$,
also denoted by $\rho_{t}$. 
The $K$-averaged coherent state transform (CST) is defined as the map 
\ba
\nonumber
C_{t} & : & L^{2}(K,dx)\mapsto {\cal H}(K_{\C})  \\
\la{211}
C_{t}(f)(g) & = & \int_{K} f(x) \rho_t (x^{-1} g)dx, \qquad  
f\in L^{2}(K,dx) , 
\ g\in K_{\C}
\ea
where $dx$ is the normalized Haar measure on $K$ and  
${\cal H}(K_{\C})$ is the 
space of holomorphic functions on $K_{\C}$.
For each $f\in L^{2}(K,dx)$, 
$C_{t}f$ is just the analytic continuation to $K_\C$ of the solution of the heat equation, 
\begin{equation}
\label{214}
\frac{1}{\pi}\frac{\partial  u}{\partial t} = \Delta u,
\end{equation}
with initial condition given by $u(0,x)=f(x)$.  
Therefore, $C_{t}(f)$ is given by 
\begin{equation}
\label{215}
C_{t}f(g) = ({\cal C}\circ \rho_t \star f)(g) = 
\left({\cal C}\circ e^{t\pi\Delta}f \right)(g),
\end{equation}
where $\star$ denotes the convolution in $K$ and ${\cal C}$ denotes analytic 
continuation from $K$ to $K_\C$. 
Let $d\nu_t$ be the $K$-averaged heat kernel measure on 
$K_{\C}$ defined in \cite{Ha1}. 
Then the following result holds.

{\theorem 
\label{t21}
{\rm [Hall]} 
For each $t> 0$, the mapping $C_t$ defined in (\ref{211}) is a unitary 
isomorphism from 
$L^{2}(K,dx)$ onto the Hilbert space $L^{2}(K_\C,d\nu_t)\bigcap \H (K_{\C})$.}
\\

To obtain a more explicit description of this CST, 
consider the expansion of $f\in L^{2}(K,dx)$
given by the Peter-Weyl theorem, 
\begin{equation}
\label{212}
f(x) = \sum_{\rm R} \mbox{tr} ({\rm R} (x)A_{{\rm R}} ),
\end{equation}
where the sum is taken over the set of (equivalence classes of) irreducible 
representations of $K$, and 
$A_{{\rm R}}\in End \ V_{\rm R}$ is given by
\begin{equation}
\label{213}
A_{{\rm R}} = \left( \dim V_{\rm R}\right) \int_K f(x) {\rm R} (x^{-1})dx, 
\end{equation}
$V_{\rm R}$ being the representation space for ${\rm R}$. 
Then one obtains:
\begin{equation}
\label{216}
C_{t}f(g) = \sum_{\rm R} e^{-t\pi c_{\rm R}}\,\mbox{tr} ({\rm R} (g)A_{{\rm R}} ),
\end{equation}
where $c_{\rm R}\geq 0$ is the eigenvalue of $-\Delta$ on functions of the type 
$$\mbox{tr} (A {\rm R}(x)), \quad A\in End(V_{\rm R}).$$


\subsection{Coherent State Transform for Class Functions on Lie groups}
\la{s21b}

From now on, let $K$ be a compact simple Lie group of $ADE$ type and let $<,>$ 
be the Ad-invariant inner product on 
$Lie(K)$ for which the roots have squared length 2. The extension of the results to 
non-$ADE$ type groups should not present a problem.
We will study here the restriction of the CST to class
(i.e. $Ad_K$-invariant) functions  and its relation to
the CST transform on a maximal  torus $T \subset K$. The main result
in this section is theorem \ref{tb1}.
Let $K/Ad_K$ be the  quotient space for the adjoint action of $K$ on itself. 
As we will show in section \ref{s42}, it turns out that the image of appropriately chosen 
distributions on $K/Ad_K \cong T/W$ 
(where $W$ is the Weyl group), 
related with Bohr-Sommerfeld conditions in geometric quantization, 
with respect to a
natural extension of $C_{t}$ in (\ref{215}) and (\ref{216}), 
gives functions satisfying 
quasi-periodicity conditions in the imaginary directions of $K_\C$.
These functions  
correspond to holomorphic sections of the pull-back 
of line bundles over certain moduli spaces of holomorphic vector bundles 
over an elliptic curve $X_\tau$. Here, the appropriate metric to be considered 
on $K_\C$ is related to the complex structure $\tau$ on $X_\tau$, where $\tau\in\C$, ${\rm Im}(\tau) >0$. 

Let $\mathfrak{h}$ be the Cartan subalgebra of $K_\C$ corresponding
to $T$, and let  $\Lambda_R, \Lambda_W \subset \h_\R^*$, $ \check \Lambda_R,
\ck \La_W \subset \h_\R$ be the root, weight, coroot and coweight 
lattices respectively. We consider fixed a choice of positive
roots. Denote also by $<,>$ the inner product induced on 
$\mathfrak{h}_\R^*$ by the inner product $<,>$ on $\mathfrak{h}$.

{}From the Peter-Weyl expansion (\ref{212}) and Schur's lemma, one sees that the space 
$\left(L^2(K,dx)\right)^{Ad_K}$ of Ad-invariant functions
on $L^2(K,dx)$ corresponds to chosing all the endomorphisms
$A_{{\rm R}}$ proportional to the identity $A_R = a_R I$. Therefore, any  
$f \in \left(L^2(K,dx)\right)^{Ad_K}$ can be expressed as 
\begin{equation}
\la{b10}
f = \sum_{\lambda \in \La_W^+} a_{\lambda} \chi_\lambda, 
\end{equation}
where we labelled irreducible representations of $K$
by the highest weights $\lambda$ in the set of 
dominant weights $\La_W^+ \subset \La_W$,
$\chi_\lambda = {\rm tr}({\rm R}_\lambda)$ is the character
corresponding to $\lambda$.

By restricting the coherent state 
transform (\ref{211}) and (\ref{216}) to the closed subspace of 
Ad-invariant functions on $K$ we obtain

{\proposition \label{p21} The restriction $C_t^{\rm Ad}$ of the CST 
(\ref{211}) to the Hilbert space 
$L^2(K,dx)^{Ad_K}$ is an isometric isomorphism
onto the Hilbert space\\ $L^{2}(K_\C,d\nu_t)^{Ad_{K_\C}}\bigcap \H (K_{\C})$.}

\begin{proof}
{}From (\ref{b10}) and (\ref{216})
we see that for $f\in (L^2(K,dx))^{Ad_K}$ we have  
\ba
\nonumber
C_t^{\rm Ad}f(g) = C_tf(g) &=&  
\sum_{\lambda \in \La_W^+} a_{\lambda}\, e^{-t\pi {c_\lambda}} \, \chi_\lambda (g)  =   \ \\ 
&=& \sum_{\lambda \in \La_W^+} a_{\lambda}\, e^{-t\pi {c_\lambda}} \,  
{\rm tr}({\rm R}_\lambda (g)), \,\,\, g\in K_\C. 
\la{b3}
\ea
where $c_\lambda = c_{{\rm R}_{\lambda}}$, and therefore the image of an $Ad_K$-invariant function on $K$
is an $Ad_{K_\C}$-invariant function on $K_\C$. 
On the other hand, a non $Ad_K$-invariant $f$ has in its Peter-Weyl expansion at least one $A_R$ which is not 
proportional to the identity and therefore its image will not be $Ad_{K_\C}$-invariant.
 Let now $F\in L^{2}(K_\C,d\nu_t)^{Ad_{K_\C}}\bigcap \H (K_{\C})$. From the ontoness of $C_t$ and from 
Schur's lemma as above, we see that $F$ has the form (\ref{b3}) and is therefore the image under $C_t$ of 
an $Ad_K$-invariant function on $K$.
\end{proof}

In view of the isomorphism of $C^\infty(K)^{Ad_K} \cong C^\infty(T)^W$ it is interesting to relate
the CST $C_t^{\rm Ad}$ for Ad-invariant functions on 
$K$ with the ``abelian'' CST $C_t^{\rm ab}$ for 
functions on $T$. Amazingly, the interesting relation 
is with $W$-anti-invariant and not $W$-invariant functions
on $T$. This fact will have important consequences
for the  application to non-abelian theta functions
that we pursue in the next sections. In particular,
it will lead to the well known shifts of level
$k \rightarrow k+h$ and of weights $\lambda \rightarrow \lambda + \rho$, where $h$ is the dual 
Coxeter number for $Lie(K)$ and $\rho\in \Lambda_W$ is the Weyl vector given by the half the sum of the 
positive roots.

First of all notice that from the Weyl integration formula (see for example \cite{Kn})
we see that there exists  an isometric isomorphism into the space of $W$-invariant functions on $T$, 
\ba
\nonumber
L^2(K,dx)^{Ad_K} & \longrightarrow & 
L^2(T,{|\sigma|^2  } dh/ |W|)^W \\
f & \mapsto & f_{|_T}
\la{b1}
\ea
where  $dh$ denotes the normalized Haar measure on $T$, $|W|$ is the order
of $W$ and $\sigma$ is the
denominator of the Weyl character formula given by
\begin{equation}
\la{b2}
\sigma (e^{2\pi i h}) = \sum_{w\in W} \epsilon(w) e^{2\pi i w(\rho)(h)}, \quad {\rm for} 
\,\,\, 
h \in {\h}_\R \ .
\end{equation}
Here, $\epsilon(w) = \det(w)$ with $w\in W$ viewed as an orthogonal transformation in ${\h}^*_\R$.  
Notice that $\sigma$ is Weyl anti-invariant, that is if 
$w\in W$ then $\sigma\circ w = \epsilon(w) \sigma$.

Let $\al_1,...,\al_l$ be a set of simple roots and $\lambda_1,...,\lambda_l$ the corresponding set of 
fundamental weights such that 
$\frac{2<\lambda_i,\alpha_j>}{<\alpha_j,\alpha_j>} = \delta_{ij}$. 
The coordinates on $\h_\R$
corresponding to the coroots ${\ck \alpha_j}$ will be denoted by $x_j$.
By considering the $T$ invariant and $W$ invariant metric 
on $T$ associated with the inner product $<,>$ we have the Laplace operator on $T$ given by 
\begin{equation}
\label{b52}
\Delta^{\rm ab} = \sum_{i,j=1}^l \frac{1}{4\pi^2} 
C^{ij} \frac{\partial^2}{\partial {x_i} \partial {x_j}}
\end{equation}
where $C^{ij} = <\lambda_i, \lambda_j>$ is the inverse Cartan matrix for $Lie(K)$.

As in the case of the group $K$, there is a unique analytic continuation $\rho_t^{\rm ab}$ to $T_\C$ of the 
heat kernel for the Laplacian $\Delta^{\rm ab}$ on $T$, and we define for each $t>0$ 
the ``abelian CST'' as the map:
\ba
C^{\rm ab}_{t} & : & L^{2}(T,dh)\rightarrow {\cal H}(T_{\C}) \nonumber \\
C^{\rm ab}_{t}(f)(z) & = & \int_{K} f(h) \rho^{\rm ab}_t (h^{-1} z)dh, \qquad  f\in L^{2}(T,dh) , 
\ z \in T_{\C}, 
\la{b4}
\ea
where ${\cal H}(T_{\C})$ is the space of holomorphic functions on $T_{\C} \cong \h/ \ck \La_R \cong 
(\C^*)^l$. It is easy to see that if $f\in L^{2}(T,dh)$ is given by 
\begin{equation}
\la{b5}
f(h) = \sum_{\lambda \in \La_W} b_{\lambda} \, e^{2\pi i \lambda}(h), \,\,\,\, {\rm for}\,\, h\in T,
\end{equation}
then
\begin{equation}
\la{b6}
C^{\rm ab}_{t}(f)(z) = 
\sum_{\lambda \in \La_W} b_{\lambda}\, e^{-t\pi{||\lambda||^2}}  
e^{2\pi i \lambda} (z),\,\,\, {\rm where} \,\, z\in T_\C.
\end{equation}

As before, Hall's result applies to this CST and we have
{\corollary 
\label{t22}
{\rm [Hall]}  For each $t> 0$, the mapping $C^{\rm ab}_t$
is an isometric isomorphism onto the Hilbert space $L^{2}(T_\C,d\nu^{\rm ab}_t)\bigcap \H (T_{\C})$,
where $d\nu^{\rm ab}_t$ is the averaged heat kernel measure on $T_\C$.
}

\vspace{2mm}

The explicit expression for $\nu^{\rm ab}_t$ is 
\begin{equation}
\label{iproclaim}
\nu^{\rm ab}_t (v) = \left( \frac{2}{t}\right)^{l/2} \sqrt{\det [C_{ij}]}\, e^{\frac{\pi}{2t}<v-\bar v,v-\bar v>},
\end{equation}
where $z=e^{2\pi i v}$, for $v \in \h$ and $T_\C \cong \h / \ck \La_R$, and conjugation $v \mapsto \bar v$ is the 
anti-linear involution in 
$\h$ preserving $\h_\R$. 

{}Functions in (\ref{b3}) are the analytic continuation
of Ad-invariant solutions of the heat equation on $K$ while
those in (\ref{b6}) are the analytic continuation
of solutions of the heat equation on $T$ and which can in turn be extended to 
$Ad_K$-invariant functions on $K$.  
As noticed by Fegan \cite{Fe} two facts
make it particularly simple to relate Ad-invariant 
solutions  of the heat equation on $K$
with Weyl {\em anti-invariant} solutions of the
heat equation on $T$.
The first is the Weyl character formula
\begin{equation}
\la{b7}
\chi_\lambda = \frac{1}{\sigma}\sum_{w \in W} \epsilon(w) e^{2\pi i w(\lambda + \rho)}
\end{equation}
and the second is the 
identity
\begin{equation}
\la{b8}
c_\lambda = {||\lambda + \rho||^2 - ||\rho||^2} \ .
\end{equation}
Let $L^{2}(T,dh)^W_-$ and 
$L^{2}(T_\C,d\nu^{\rm ab}_t)^W_-$ denote the
$W$-anti-invariant subspaces. Since the actions of $W$ and $\Delta$ on $L^2(T,dh)$
commute,  we have a result analogous to Proposition
\ref{p21} (with similar proof which we omit):

{\proposition \label{p22} The restriction  of the abelian CST 
(\ref{b4}) (which we will denote by the same symbol
$C_t^{\rm ab}$) to the space 
$L^{2}(T,dh)^W_-$
is an isometric isomorphism
onto the Hilbert space $L^{2}(T_\C,d\nu^{\rm ab}_t)^W_- \bigcap \H (T_{\C}) $.}

\vspace{2mm}

{}From (\ref{b3}), (\ref{b6}), (\ref{b7}) and (\ref{b8}) we see
that by multiplying an Ad-invariant solution 
of the heat equation on $K$  by $\sigma e^{- t \pi ||\rho||^2}$
we obtain a $W$-anti-invariant solution of the
heat equation on $T$. Moreover, this map takes
$C_t^{\rm Ad}$ to $C_t^{\rm ab}$ restricted
to $W$ anti-invariant functions.
More precisely, consider the maps 
$\varphi$ and $\varphi_\C$  given respectively by
\begin{equation}
\begin{array}{cccc}
\nonumber
\varphi  : & L^2(K,dx)^{Ad_K} & \rightarrow & L^{2}(T,dh)^W_- \\[1mm]
&  f & \mapsto &  \frac{\sigma}{\sqrt{|W|}}    f_{|_T}
\la{b12}
\end{array}
\end{equation}
and
\begin{equation}
\begin{array}{cccc}
\nonumber
\varphi_\C \ : & L^{2}(K_\C,\nu_t)^{Ad_{K_\C}}\bigcap \H (K_{\C}) & \rightarrow &
L^{2}(T_\C,d\nu^{\rm ab}_t)^W_- \bigcap \H (T_{\C})    \\[2mm]
& \ f &\mapsto& e^{- t \pi {||\rho||^2}} \frac{\sigma}{\sqrt{|W|}}  \ f_{|_{T_\C}}
\la{b11}
\end{array}
\end{equation}
where $\sigma$  denotes also the analytic continuation of (\ref{b2}). 
We then have the following

\begin{theorem}
\la{tb1} 
The maps $\varphi$ and $\varphi_\C$ are
isometric isomorphisms and
the following diagram is commutative
\begin{equation}
\begin{array}{ccc}
L^2(K,dx)^{Ad_K} &  \stackrel{C_t^{\rm Ad}}{\rightarrow } &
L^{2}(K_\C,d\nu_t)^{Ad_{K_\C}}\bigcap \H (K_{\C}) \\
\downarrow \varphi & & \downarrow \varphi_\C \\
 L^{2}(T,dh)^W_- &  \stackrel{C_t^{\rm ab}}{\rightarrow } &
L^{2}(T_\C,d\nu^{\rm ab}_t)^W_- \bigcap \H (T_{\C}) 
\end{array}
\la{b9}
\end{equation}
\end{theorem}

\begin{proof}
{}From the Weyl integration formula we see that $\varphi$ is an isometry. 
It is easy to check that if $f\in L^{2}(T,dh)^W_-$, then its expansion 
(\ref{b5}) receives contributions only from non-singular weights $\lambda \in \Lambda_W$, 
i.e. such that \newline $<\lambda,\alpha_i>\neq 0$ for all simple roots $\alpha_i$.
This implies that such $f$ is of the form
$$
f(h)=\sum_{\stackrel{\lambda \in \Lambda_W}{\lambda{\rm\  non-singular}}} 
b_\lambda e^{2\pi i \lambda(h)}.
$$
Using Weyl anti-invariance and the fact that any regular $\lambda'\in\Lambda_W^+$
is of the form $\lambda'=\lambda+\rho, \lambda \in \Lambda_W^+$, we obtain
$$
f(h)= \sum_{w\in W} \epsilon (w)\sum_{\lambda \in \Lambda_W^+} b_{\lambda + \rho} e^{2\pi i w(\lambda+\rho)(h)} =
\sigma(h) \sum_{\lambda \in \Lambda_W^+} b_{\lambda + \rho} \chi_\lambda (h), 
$$
so that $f/\sigma$ can be extended to $L^2(K,dx)^{Ad_K}$, and $\varphi$ is an isomorphism. 
On the other hand, 
from (\ref{b3}), (\ref{b6}), (\ref{b7}) and (\ref{b8})
we see that the diagram commutes. It then follows from Propositions 
\ref{p21} and \ref{p22} 
that $\varphi_\C$ is also an isometric isomorphism.
\end{proof}


\subsection{Extension to Distributions}
\label{vemdo4.2}

In order to apply later the CST to the study of non-abelian theta functions 
on an elliptic curve with modular parameter $\tau$ in the Teichmuller space of genus 1 curves 
${\cal T}_1 = \{\tau \in \C: {\rm Im}\tau >0\}$, let us 
consider the complex non-self-adjoint Laplacian on $K$
\begin{equation} 
\la{422}
\Delta^{(-i\tau)} = -i \tau \Delta = -i \tau \sum_{j=1}^{\dim K} X_jX_j \ . 
\end{equation}
where 
$\{X_i, \, i=1,...,\dim K\}$ is an  orthonormal 
basis for $Lie(K)$ viewed as the space 
of left-invariant vector fields on $K$ as in section \ref{s21a}.
Then we have
{\proposition \la{p41}
For each $\tau \in {\cal T}_1$ and each $t> 0$, the mapping $C_t^{\tau}$
\begin{equation}
\la{423}
C_t^{\tau} = {\cal C} \circ e^{{t \pi}{\Delta^{(-i\tau)} }} \ 
: \ L^2(K, dx) \rightarrow \H (K_\C)
\end{equation}
is a unitary isomorphism onto the Hilbert 
space 
$$L^{2}(K_\C,d\nu_{t\tau_2})\bigcap \H (K_\C),$$
where $d\nu_{t\tau_2}$ is the averaged heat kernel measure
corresponding  to the Laplacian $\Delta$ on $K$.
The restriction of the CST   
(\ref{423}) to the space 
$L^2(K,dx)^{Ad_{K}}$, also denoted by $C_t^{\tau}$, is an isometric isomorphism
onto the Hilbert space
$$L^{2}(K_\C,d\nu_{t\tau_2})^{Ad_{K_\C}}\bigcap \H (K_\C).$$}
\begin{proof}
Let $\tau_1 = {\rm Re}(\tau)$ and 
decompose the transform (\ref{423}) as 
$$
{\cal C} \circ e^{{t \pi}{\Delta^{(\tau_2)} }} 
\circ e^{{t \pi}{\Delta^{(-i \tau_1)} }} = C_t^{i\tau_2} \circ e^{{-i\tau_1 t \pi}
{ \Delta}}. 
$$ 
The Laplace operator $\Delta$ is self-adjoint on $L^2 (K,dx)$ and therefore the 
operator 
$$
e^{-i\tau_1 t\pi\Delta}:L^2 (K,dx) \longrightarrow L^2 (K,dx) 
$$
is unitary. The unitarity of $C_t^{i\tau_2}$ follows from Theorem \ref{t21}. 
To obtain the second statement notice 
that the image under the CST (\ref{423}) of $f$ in the form (\ref{b10}) is given by
\begin{equation}
\label{424}
C_t^{\tau}f(z) =   
\sum_{\lambda \in \La_W^+} a_{\lambda}\, e^{i\pi\tau t {c_\lambda}}  
\chi_\lambda (z). 
\end{equation}
The proof then follows from the proof of proposition \ref{p21} with obvious changes.
\end{proof}

We will now extend the CST transform of proposition \ref{p41} to the space of 
distributions  $C^{\infty}(K)^\prime$ \cite{V}. 
This is the space of Fourier series of the form \cite{Sc}
\begin{equation}
\label{425}
f = \sum_{\lambda\in\La_W^+} {\rm tr} ({\rm R}_\lambda A_\lambda)
\end{equation}
for which there exists an integer $N>0$ such that 
\begin{equation}
\label{429}
 \lim_{||\lambda||\rightarrow\infty} 
\frac{||A_\lambda||}{(1+||\lambda||^2)^N} = 0,
\end{equation}
where $||A_\lambda||$ is the operator norm of the endomorphism 
$A_\lambda\in End({\rm R}_\lambda)$. 
Consider the extension of the adjoint action of $K$ on
$C^{\infty}(K)$ to $C^{\infty}(K)^\prime$ defined by 
$$
x \cdot f = \sum_{\lambda\in\La_W^+} {\rm tr} ({\rm R}_\lambda A^x_\lambda) \ , 
$$
where $A^x_\lambda = R_\lambda(x) A_\lambda R_\lambda(x^{-1}), \forall x\in K$.
If $f$ is in the space of 
$Ad_{K}$-invariant distributions 
$f\in (C^{\infty}(K)^\prime)^{Ad_{K}}$ then Schur's lemma and (\ref{425}) again imply 
that it has a unique representation of the form
\begin{equation}
\label{427}
f = \sum_{\lambda\in\La_W^+} a_\lambda \chi_\lambda .
\end{equation}

The Laplace operator and its powers act as continuous linear operators on
the space  $C^{\infty}(K)^\prime$ and for $\tau_2 > 0$ define the 
action of the operator $e^{\pi t \Delta^{(-i\tau)}}$ on it. 
For $f$ of the form (\ref{425}) we have
\begin{equation}
\label{426}
e^{\pi t \Delta^{(-i\tau)}} f  = \sum_{\lambda\in\Lambda_W^+} 
e^{i \pi t \tau c_\lambda} {\rm tr} ({\rm R}_\lambda A_\lambda).
\end{equation} 

\begin{proposition}
\label{p42}
If $f = \sum_{\lambda\in\La_W^+} {\rm tr} ({\rm R}_\lambda A_\lambda) 
\in C^{\infty}(K)^\prime$ 
then the series 
\begin{equation}
\label{426.5}
\sum_{\lambda\in\La_W^+} 
e^{i \pi t \tau c_\lambda} {\rm tr} ({\rm R}_\lambda (g)A_\lambda)
\end{equation}
where $g\in K_\C$, defines a holomorphic function on $K_\C \times {\cal T}_1$ which we denote by 
$\left( {\cal C}\circ e^{\pi t \Delta^{(-i\tau)}}\right)(f)$.
\end{proposition}

\begin{proof}
The growing condition (\ref{429}), together with (\ref{b8}), implies that there exists a $c>0$ such that 
$$
||A_\lambda ||e^{-t\pi \tau_2 c_\lambda} \leq e^{-c ||\lambda||^2},
$$  
for $||\lambda ||$ sufficiently large. 
On the other hand, writing $g = x \exp(iY)$ for $x \in K, Y\in Lie(K)$, one has 
$||R_\lambda (g)|| \leq \exp (M ||Y|| ||\lambda||)$, for some constant $M>0$.
Therefore, the series (\ref{426.5}) 
is uniformly convergent on compact subsets 
of $K_\C \times {\cal T}_1$ and its sum defines a holomorphic function there.
\end{proof}

\begin{corollary} \la{co41}
 If the distribution $f$ is $Ad_{K}$-invariant
 of the form (\ref{427}) then the series 
\begin{equation}
\label{428}
\sum_{\lambda\in\La_W^+} a_\lambda e^{i \pi t \tau c_\lambda} \chi_\lambda (g)
\end{equation}
defines an $Ad_{K_\C}$-invariant holomorphic function on  $K_\C \times {\cal T}_1$.
\end{corollary}

\begin{definition}
\label{def41}
The $K$-coherent state transform for the elliptic curve 
$X_\tau= \C/(\Z\oplus \tau\Z)$, $\tau \in {\cal T}_1$, and $t>0$ is the map
\begin{equation}
\label{428.5}
C_t^{\tau}={\cal C} \circ e^{\pi t \Delta^{(-i\tau)}} \ : 
\left(C^{\infty}(K\right)^\prime)^{Ad_{K}} 
\longrightarrow {\cal H}(K_\C)^{Ad_{K_\C}}.
\end{equation} 
\end{definition}

\begin{remark}
The role of the elliptic curve $X_\tau$ will become clear in section \ref{s4}.
\end{remark}

We will also need the extension of the abelian coherent state transform in (\ref{b4}) 
to the present case.
This will be related later on in definition \ref{d23.1} to a
CST associated to the matrix $\widetilde 
\Omega = \tau C^{-1}$ and defined by 
$$
C^{{\rm ab}(\tau)}_t = {\cal C} \circ  e^{-i\pi \tau \Delta^{{\rm ab}}}: \,\,C^{\infty}(T)' \rightarrow 
{\cal H} (T_\C).
$$

Consider the maps $\varphi$ and $\varphi_\C$  generalizing the maps 
in (\ref{b12}) given by
\ba
\nonumber
\left(\C^{\infty}(K)'\right)^{Ad_{K}} & \rightarrow 
&\left(\C^{\infty}(T)'\right)^W_- \\
\varphi \ : \ f=\sum_{\lambda \in \Lambda_W^+}a_\lambda \chi_\lambda  &\mapsto&  
\frac{1}{\sqrt{|W|}} \sum_{w \in W}\epsilon(w) \sum_{\lambda \in \Lambda_W^+} 
a_{\lambda} e^{2\pi i w(\lambda + \rho)}
 \la{car1}
\ea
and
\ba
\nonumber
\H (K_\C) & \rightarrow & \H (T_{\C})^W_-    \\
\varphi_\C \ : \ f &\mapsto& e^{i t \tau \pi {||\rho||^2}}
\frac{\sigma}{\sqrt{|W|}}  \ f_{|_{T_\C}}.
\la{car2}
\ea
We then have the following

\begin{theorem} 
\la{pcar} 
The maps $\varphi$ and $\varphi_\C$ are isomorphisms,  
and the following diagram is commutative
\begin{equation}
\begin{array}{ccc}
\left(\C^{\infty}(K)'\right)^{Ad_{K}} &  \stackrel{C_t^{\tau}}{\rightarrow } &
{\cal H}(K_\C)^{Ad_{K_\C}} \\
\downarrow \varphi & & \downarrow \varphi_\C \\
\left(\C^{\infty}(T)'\right)^W_- &  \stackrel{C_t^{\rm ab(\tau)}} {\rightarrow } &
\H (T_{\C})^W_- 
\end{array}
\la{car3}
\end{equation}
and the maps $C_t^{\tau}$ and $C_t^{\rm ab(\tau)}$ are injective.
\end{theorem}

\begin{proof}
A similar argument to the one in the proof of theorem (\ref{tb1}),
the Weyl character formula (\ref{b7}) and (\ref{427}) show that $\varphi$ is an isomorphism. 
On the other hand, from (\ref{424}), (\ref{427}) and (\ref{b6}) we see that 
$C_t^{\tau}$ and $C_t^{\rm ab(\tau)}$ 
are injective. Indeed, two distributions $f_1$ and $f_2$ with representations 
(\ref{427}) are different if and only if
there exists $\lambda\in \La_W^+$ such that the corresponding coefficients
$a_\lambda^1$, $a_\lambda^2$ are different.
Then, the two holomorphic functions $C_t^{\tau}f_1$ and $C_t^{\tau}f_2$ 
have also different coefficients with respect to 
$\chi_\lambda$ and are therefore different, as can be readily seen by restriction to $K$. 
The injectivity of $C_t^{{\rm ab}(\tau)}$ is proved in a similar way.
The fact that $\varphi_\C$ is also an isomorphism, follows from lemma 9 of \cite{Ha1}.
Indeed, any $f\in {\cal H}(K_\C)$ has a unique Peter-Weyl decomposition which
corresponds to the unique analytic continuation to $K_\C$ of $f|_K\in L^2(K,dx)$.
(The same applies to $T_\C$.)
Therefore, again, an analogous argument to the one in the proof of theorem (\ref{tb1})
shows that $\varphi_\C$ is an isomorphism.
Finally, the diagram commutes due to the Weyl character formula, (\ref{b8}) and the 
expressions for the CST in (\ref{424}), (\ref{b6}) and 
(\ref{d23.1}).
\end{proof}



\section{Coherent State Transform and Theta Functions on 
Polarized Abelian Varieties}
\label{s22}

In this section, we extend the results of \cite{FMN} to abelian varieties with general polarization.
This extension, together with the results of the previous section, will be applied in the next sections 
to the study of non-abelian theta functions in genus one. 

Let $V$ be an $l$-dimensional complex vector
space and $\Lambda \cong \Z^{2l}$ a maximal lattice in $V$
such that the quotient
\begin{equation}
\label{eq2.2.1}
M = V / \Lambda
\end{equation} 
is an abelian variety, i.e. a complex torus which
can be holomorphically embedded in projective space. 
For later convenience we will assume 
that $M$ is endowed with a polarization $H$, not necessarilly principal 
\cite{GH, BL, Ke}. 
By
definition, $H$ is a positive definite hermitean form on $V$ whose imaginary
part is integral on $\Lambda .$ Let $E=-{\rm Im}H,$ so that $E$ is an
integral alternating bilinear form on $\Lambda .$ According to the
elementary divisor theorem, there is a canonical basis of $\Lambda ,$ $\beta
_{1},\ldots ,\beta _{l},\tilde{\beta}_{1},\ldots ,\tilde{\beta}_{l}$
characterized by 
\begin{eqnarray}
\label{2.2.2} 
\nonumber 
E(\beta _{i},\beta _{j}) &=&E(\tilde{\beta}_{i},\tilde{\beta}_{j})=0 \\
E(\beta _{i},\tilde{\beta}_{j}) &=&\delta_{i}\delta _{ij},\qquad i,j=1,\ldots ,l,
\end{eqnarray}
where $\delta_1|\delta_2 | ... |\delta_l$ are positive integers, depending only on 
$E,$ and $\delta _{ij}$ is Kronecker's delta symbol. Define $E_1 = - {\rm Im}H_1$ to be 
the form with $\delta_1 =1$.

Now, let us decompose $V$ into isotropic subspaces with respect to $E_1,$ 
$V_{1}=\oplus _{i=1}^{l}\Bbb{R}\beta _{i},$ $V_{2}=\oplus _{i=1}^{l}\Bbb{R}
\tilde{\beta}_{i},$ and decompose the lattice in the same way $\Lambda
_{i}=V_{i}\cap \Lambda ,$ $i=1,2,$ so that $\Lambda =\Lambda _{1}\oplus
\Lambda _{2}.$ Let $\alpha $ be the semicharacter for $H$ which is trivial
on $\Lambda _{1}$ and $\Lambda _{2},$ i.e, $\alpha $ is the unique map 
$\alpha :\Lambda \rightarrow U(1)$ satisfying 
\[
\left\{ 
\begin{array}{c}
\alpha (\lambda +\lambda ^{\prime })=\alpha (\lambda )\alpha (\lambda
^{\prime })e^{-\pi iE_1(\lambda ,\lambda ^{\prime })} \\ 
\alpha |_{\Lambda _{1}}=\alpha |_{\Lambda _{2}}=1.
\end{array}
\right. 
\]
To this particular Appell-Humbert pair $(\alpha ,H_1)$, we can associate a
line bundle which we will denote by $L_{1}=L(\alpha,H_1)$ over $M$ via the
following (canonical) factors of automorphy 
\[
a(v,\lambda)=\alpha (\lambda )e^{\pi H_1(v,\lambda)+\frac{\pi }{2}H_1(\lambda
,\lambda )}.
\]
Recall that, via the canonical identification of $H^2(M,\Z)$ with the space of
integral alternating bilinear forms on $\Lambda$, the first Chern class
of $L_1$, $c_1(L_1)$, corresponds to the form $E_1$.

Level $k$ theta functions on $M$ are holomorphic sections, $\tilde \theta$, of 
$L_{1}^{k}= L (\alpha^k, kH_1 )$, $\tilde \theta \in H^0(M,L_1^k)$. 

For convenience, we will consider the
following different but equivalent factors of automorphy for $L_{1}.$ Let $S$
be the $\Bbb{C}$-bilinear extension of $H|_{V_{1}\times V_{1}}$ to $V\cong 
\Bbb{C}\cdot V_{1}$, and define 
$$
F = \frac{1}{2i}(H-S)
$$ 
which is a form $\C$-linear in the first variable. 
Then, $L_{1}^k$ is also given by the following (classical) factors of automorphy
$e(\lambda,v)=\alpha(\lambda)e^{2\pi ik F(v,\lambda)+ \pi i kF(\lambda,\lambda)}$, 
which can be rewritten as
\begin{equation}
\label{2.2.6}
e(\lambda _{1}+\lambda _{2},v)=e^{2\pi ikF(v,\lambda _{2})+\pi ikF(\lambda
_{2},\lambda _{2})},\qquad \lambda _{i}\in \Lambda _{i}.
\end{equation}

To write down explicit expressions for the theta functions,   
consider the lattice dual to $\Lambda$ with respect to $E_1$,  
\begin{equation}
\label{2.2.3}
\widehat{\Lambda} = \{ v\in V: E_1 (v, \lambda)\in \Z, \, 
\forall \lambda\in\Lambda\} = \widehat{\La_1}\oplus \widehat{\La_2} \supset\Lambda
\end{equation}
where $\widehat{\La_i} = V_i \bigcap \widehat{\Lambda}\supset \La_i$, for $i=1,2$. 
A basis for $\widehat{\La_1}$ and $\widehat{\La_2}$ is given respectively by 
\begin{eqnarray}
\nonumber
\b^{'}_{j} & = & \frac{1}{\delta_j} \b_j,  \\
\nonumber
\tilde \b^{'}_{j} & = & \frac{1}{\delta_j} \tilde \b_{j} = 
\sum_{i = 1}^l \Omega_{ij} \b_i, \quad 
j = 1, ..., l,
\label{2.2.4}
\end{eqnarray}
where $\Omega = (\Omega_{ij})$ is a matrix in the Siegel upper 
half space $\HH_l$ of 
symmetric $l\times l$ matrices with positive imaginary part. 

One computes,  
\begin{equation}
\label{2.2.7}
F(\tilde\b'_{i}, \tilde\b'_{j}) = -\Omega_{ij}.
\end{equation}

Using the automorphy factors (\ref{2.2.6}), we see that the space 
\break $H^0(M,L_1^k)$, is isomorphic to the space of holomorphic 
functions on $V/\La_1 \cong (\C^*)^l$ satisfying quasi-periodicity 
conditions in the directions of $\La_2$ given by 
\begin{equation} 
\label{2.2.8}
\theta (v + b) = e^{2\pi i kF(v,b) + \pi i kF(b,b)} \theta (v), \,b\in \La_2.
\end{equation}

Let us denote the latter space by ${\cal H}_{k,\Omega}$, 
$$
{\cal H}_{k,\Omega}\subset {\cal H}((\C^*)^l),
$$
where ${\cal H}((\C^*)^l)$ denotes the space of holomorphic functions 
on $(\C^*)^l$.
Conditions (\ref{2.2.8}) then imply that there is one independent 
theta function $\theta\in {\cal H}_{\Omega,k}$, 
for every $C\in  (k^{-1} {\widehat \La_2})/\La_2$ given by 
\begin{eqnarray}
\nonumber
\theta_C (v) & = & \sum_{b\in \La_2} e^{-\pi i kF(c_0 +b, c_0 +b) - 
2\pi i kF(c_0 +b, v)}\\
 & = & \sum_{c\in C} e^{-\pi i kF(c,c) - 2\pi i kF(c, v)}, 
\label{2.2.9}
\end{eqnarray}  
where $[c_0] = C \in (k^{-1} {\widehat \La_2})/\La_2$. In coordinates, 
\begin{equation}
\la{330}
\begin{array}{lll}
c  =  \sum_{j} m_j \frac{\tilde\b'_{j}}{k}, && 
m=(m_1,...,m_l)\in \Z^l , \\
b  =  \sum_j p_j \tilde\b_{j} = \sum_j p_j \delta_j 
\tilde\b'_{j}, &&  p=(p_1,...,p_l)\in\Z^l ,\\
v  =  \sum_j z_j \b_j, &&  z= (z_1,...,z_l)\in \C^l ,
\end{array}
\end{equation}
the theta functions take the form
\begin{equation}
\theta_m (z,\Omega) = \sum_{p\in \Z^l} e^{\pi i (m+k\delta p)\cdot 
\frac{\Omega}{k}  (m+k\delta p) + 
2\pi i (m+k\delta p)\cdot z}, 
\label{2.2.10}
\end{equation}
where $m\in \Z^l/k(\delta_1 \Z\oplus\cdots\oplus\delta_l\Z)$, with 
$\delta_1 = 1$, $\delta = {\rm diag}(\delta_1, ...,\delta_l)$ 
and $z\cdot z' = z_1z'_1+\cdots + z_lz'_l$.
In these coordinates the automorphy factors (\ref{2.2.6})
take the form
\begin{equation}
\la{2.32a}
e(\tilde\b_{j},z)=e^{-2\pi i k z_j \delta_j - \pi i k \Omega_{jj}\delta_j^2}.
\end{equation} 

For the benefit of section \ref{s4}  let 
us consider a more general basis $\{\gamma_j\}_{j=1}^{2l}$ (not necessarily canonical) 
for $\La_1\oplus \La_2$ and
its dual basis $\{\gamma^\prime_j\}_{j=1}^{2l}$ for $\wh \La_1\oplus \wh \La_2$. 
Notice that 
the basis of $\La_1$ given by $\{\gamma_j\}_{j=1}^{l}$ can be extended to a 
canonical basis 
of $\La_1\oplus \La_2$ if and only if $\delta_j \gamma_{j+l}'$ is a basis of $\La_2$.
Let $\{ \beta_j, \tilde\beta_{j} \}$ 
continue to denote a canonical basis of $\La_1\oplus \La_2$ and 
\begin{equation}
\begin{array}{rclrcl}
\la{331}
\beta_j &=& \sum_{i=1}^l \gamma_i P_{ij}  \quad  &  \tilde\beta_{j} &=& \tilde\beta'_{j} 
\delta_{j} \\
&&&&&\\
\beta_j &=&  \b'_j \delta_{j}  \quad  & \tilde\beta'_{j} &=& \sum_{i=1}^l \gamma'_{i+l} 
Q_{ij} \\
&&&&&\\
\gamma_{j+l} &=& \sum_{i=1}^l \gamma'_{i+l} R_{ij} &&&  
\end{array}
\end{equation}
where $P,Q\in SL(l,\Z)$ and duality demands that $P^tQ =Id$. We then have 
\begin{equation}
\la{332}
\gamma'_{j+l} = \sum_{i=1}^l \gamma_i \widetilde \Omega_{ij}     
\end{equation}
where 
\begin{equation}
\la{333}
\widetilde\Omega = P \Omega P^t \in \HH_l,
\end{equation}
and also $F(\gamma_j,\gamma'_{i+l})=-\delta_{ji}$ and 
$F(\gamma'_{j+l},\gamma'_{i+l})=-\wt\Omega_{ji}$, where $\delta_{ji}$.
Considering then instead of (\ref{330}) the coordinates 
\begin{equation}
\la{334}
\begin{array}{lll}
c  =  \sum_{j} \tilde m_j \frac{\gamma'_{j + l}}{k}, && 
\tilde m=(\tilde m_1,...,\tilde m_l)\in \Z^l , \\
b  =  \sum_j \tilde p_j \gamma_{j + l} = \sum_j \gamma'_{j + l}R_{ij} 
\tilde p_i, &&  \tilde p=(\tilde p_1,...,\tilde p_l)\in\Z^l ,\\
v  =  \sum_j \tilde z_j \gamma_j, &&  \tilde z= (\tilde z_1,...,\tilde z_l)\in \C^l ,
\end{array}
\end{equation}
we obtain for the theta functions the expressions 
\begin{equation}
\la{335}
\theta_{\tilde m} (\tilde z, \wt \Omega) =\theta_{m} (z, \Omega) = 
\sum_{\tilde p \in \Z^l} 
e^{\pi i (\tilde m + k R \tilde p)\cdot \frac{\wt\Omega}{k}(\tilde m + k R \tilde p) + 
2\pi i (\tilde m + k R\tilde p)\cdot \tilde z }.
\end{equation} 

In these coordinates the automorphy factors take the form
\begin{equation}
\la{336}
e(\gamma_{l+j},\tilde z) = e^{2\pi i k (R^t \tilde z)_j + 
\pi i k (R^t \tilde\Omega R)_{jj}}.
\end{equation}

Returning to the canonical coordinates we see from (\ref{2.2.10}) 
that the theta functions are the analytic continuation
to $(\C^*)^l \cong V/\La_1$ of solutions (with $t=1/k$) of the heat equation on $U(1)^l$
\begin{equation}
\la{2.321}
\frac{1}{\pi} \frac{\partial u}{\partial t} = \Delta^{{\rm ab}(-i\Omega)} u \ ,
\end{equation}
where
\begin{equation}
\la{2.322}
\Delta^{{\rm ab}(-i\Omega)} = - \sum_{j,j'}^l \frac{i}{4\pi^2}\Omega_{jj'}
\frac{\partial^2}{\partial x_j \partial x_{j'}} \ ,
\end{equation}
and the $x_i\in [0,1]$ are angular coordinates on $U(1)^l$.
Of course, that in the coordinates (\ref{334}) the functions $\theta_{\tilde m}$ are 
the analytic continuations of solutions of 
the heat equation of (\ref{2.322}) with  $\wt\Omega$ in the place of $\Omega$.

Following \cite{FMN} we extend the CST to  distributions on $(S^1)^l$.

\begin{definition}
\la{d23.1}
The CST for the  matrix $\Omega \in \HH_l$
is the transform
\ba
\la{2.323}
C_t^{{\rm ab}\,\Omega} \ : \left( C^\infty((S^1)^l)\right)^\prime & 
\rightarrow & \H((\C^*)^l)\\
\nn
f & \mapsto & {\cal C} \circ e^{\pi t \Delta^{{\rm ab}(-i\Omega)}} f \ .
\ea
\end{definition} 
The fact that the map (\ref{2.323}) is well defined follows from Lemma 4.1
of \cite{FMN}.
The averaged heat kernel function $\nu_t^{\rm ab}$ in this case reads
\begin{equation}
\la{3.324}
\nu_t^{\rm ab}(z) = \left(\frac{2}{t}\right)^{\frac{l}{2}}\left( {\rm det} W \right)^{\frac{1}{2}}
e^{\frac{\pi}{2t} \sum_{ij}(z_i - \overline z_i) W_{ij}(z_j - \overline z_j)} \ ,
\end{equation}
where $W = \Omega_2^{-1} = {\rm Im}(\Omega)^{-1}$. In periodic coordinates $(\eta , \xi)$ dual to the 
canonical basis $\{\beta_j, \tilde\beta_{j}\}$ and related to $z$ in (\ref{330}) by  
$z = \eta + \Omega \delta \xi$, the heat kernel measure reads
\begin{equation}
\la{3.325}
\nu_t^{\rm ab}(\eta , \xi) d\eta d\xi= \left(\frac{2}{t}\right)^{\frac{l}{2}}
\left( {\rm det} \Omega_2 \right)^{\frac{1}{2}} (\det \delta)
e^{-{2 \pi}{t} \sum_{ij}\xi_i \delta_i \Omega_{2_{ij}}\xi_j \delta_j} d\eta d\xi.
\end{equation}

Consider the distributions on $(S^1)^l$ given by
\begin{equation}
\la{3.326}
\theta_m^0(x) = \theta_m(x,0) = \sum_{p\in \Z^l} e^{2\pi i (m + k\delta p)\cdot x } \ ,
\end{equation}
and let ${\cal I}_{k,\delta}$ denote the $(\delta_1 \cdots \delta_l k^l)$-dimensional
subspace of $ \left( C^\infty((S^1)^l)\right)^\prime$ generated by 
these distributions with inner product $(\cdot , \cdot)$ for which the distributions
in (\ref{3.326}) form an orhonormal basis. The distributions in ${\cal I}_{k,\delta}$ are linear combinations of Dirac delta 
distributions supported on points arising from Bohr-Sommerfeld conditions \cite{Sn,Ty}. In fact, 
\begin{equation}
\la{yeye}
\delta_{m'}(x) = \delta (x-\frac{\delta^{-1}m'}{k}) = \sum_{0\leq m_\alpha < k\delta_\al} 
e^{-2\pi i m'\cdot\frac{\delta^{-1}}{k}m} \theta^0_m (x).
\end{equation}

The extension of the results of \cite{FMN} to general polarizations
can be summarized in 

\begin{theorem}
\la{t23.1}

\noindent
\begin{itemize}

\item[1.] The image under the CST $C_{t=1/k}^{{\rm ab}\,\Omega}$ of ${\cal I}_{k,\delta}$
gives the space $\H_{k,\Omega}$ of all level $k$ theta functions
on the abelian variety $M$ with polarization given by $\delta$.

\item[2.] The function $\nu_t^{\rm ab}$ is, for $t=1/k$, the pull-back
from $M$ to $(\C^*)^l$ of a hermitean structure on $L_1^k$.

\item[3.] Consider on $\H_{k,\Omega}$ the inner product induced by the CST transform
\begin{equation}
\la{3.327}
<\theta,\theta'> = \int_{[0,1]^l \times [0,1]^l} \overline \theta \theta' \nu_{1/k}^{\rm ab}
(\eta,\xi) d\eta d\xi
\ .
\end{equation}
The CST transform of definition \ref{d23.1} is, for $t=1/k$, a unitary transform 
between $\left(\I_{k,\delta}, (\cdot , \cdot ) \right)$ and  $\left(\H_{k,\Omega}, <\cdot , \cdot > \right)$.

\end{itemize}

\end{theorem}

\begin{proof}
\noindent
\newline
1. This follows immediately from the definition of ${\cal I}_{k,\delta}$.
\newline
2. This means that  $\nu_{\frac{1}{k}}^{\rm ab}$ is the pull-back of a section of $(L_1^k \otimes {\bar L_1^k})^*$ and so 
should satisfy the quasi-periodicity conditions
\begin{equation}
\la{46p}
\nu_{1/k}^{\rm ab} (z + \tilde \b_{j}) = |e(\tilde\b_{j},z)|^{-2} \nu_{1/k}^{\rm ab} (z) =
|e^{-2\pi i k z_j \delta_j - \pi i k \Omega_{jj}\delta_j^2}|^{-2} \nu_{1/k}^{\rm ab} (z).
\end{equation}
This is easy to check by taking $\xi_j \rightarrow \xi_j +1$ in the expression (\ref{3.325}). 
\newline
3. This follows from a direct computation:
\ba
\nonumber
&&\int_{[0,1]^l \times [0,1]^l} {\bar \theta_m} \theta_{m'} \nu^{\rm ab}_{1/k}(\eta,\xi) 
 d\eta d\xi  = \delta_{mm'} (2k)^{l/2} (\det \Omega_2)^{1/2} {\rm det} (\delta)\times\\
\la{46pp} 
&& \times \sum_{p\in \Z^l}\int_{[0,1]^l} 
e^{-2\pi  k(\delta\xi +\frac{1}{k}(m+k\delta p)) \cdot \Omega_2  (\delta\xi +\frac{1}{k}(m+k\delta p))}d\xi 
= \delta_{mm'},
\ea
for all $m,m' \in \Z^l/k(\delta_1 \Z\oplus\cdots\oplus\delta_l\Z)$.
\end{proof}

We end this section by noting that the $(\delta_1\cdots\delta_l k^l)\times(\delta_1\cdots\delta_l k^l)$ matrix
$$
A_{m'm} = \left( e^{-2\pi i m' \cdot \frac{\delta^{-1}}{k}m}\right)
$$  
in (\ref{yeye}) satisties $\bar A^t A = \delta_1 \cdots \delta_l k^l I$, where $I$ is the 
$(\delta_1\cdots\delta_l k^l)$-dimensional identity matrix. This means that the distributions 
$(1/\sqrt{\delta_1 \cdots \delta_l k^l})\delta_m (x)$ in (\ref{yeye}) are orthonormal in 
$\left(\I_{k,\delta}, (\cdot , \cdot ) \right)$.


\def\G{GL(n,\C)}

\section{Vector Bundles on Elliptic Curves}

\label{s3}

In this section we recall the existence of a moduli space $\mathcal{M}_{n}=%
\mathcal{M}_{n}\left( X_{\tau }\right) $, parametrizing $S$-equivalence
classes of semistable bundles of rank $n$ and trivial determinant over an
elliptic curve $X_{\tau }$, $\tau \in {\cal T}_1$. Non-abelian theta functions of genus one are
then defined to be holomorphic sections of the line bundles over $\mathcal{M}_{n}$. 
We then define the Schottky map associated to a general
complex linear group $G$ and a Riemann surface $X$, compute it explicitly 
for the case $G=\Bbb{C}^{*},$ and relate it
to a map used in \cite{FMN} to study abelian theta functions. 
We compute the Schottky map for the group $SL(n,{\mathbb C})$
over the elliptic curve $X_{\tau }$ and use it in section \ref{s4} to pull-back
to $SL(n,{\mathbb C})$ sections of line bundles over $\mathcal{M}_{n}$ and
compare them to the CST considered before. This will relate the CST for the
elliptic curve $X_{\tau }$, in definition \ref{def41}, to non-abelian theta functions of genus one.


\subsection{The Moduli Space of Semistable Vector Bundles}

\label{s3.1}Let $X$ be a smooth complex projective algebraic curve of genus $
g$. Here we will consider \emph{holomorphic vector bundles }over $X,$ which
we will denote simply by the term\emph{\ bundle.} For details of the
constructions and proofs of the results in this section, we refer to \cite{A}, \cite{Tu}.

To construct a moduli space for bundles over $X$ one introduces the following
notions.
A bundle $E$ is called \emph{stable} (resp. \emph{semistable}) if for every
proper subbundle $F\subset E$ we have 
$\mu _{F}<\mu _{E}\quad (\mbox{resp.}\quad \mu _{F}\leq \mu _{E}),$ 
where $\mu _{E}$ denotes the \emph{slope} of a bundle $E$, defined by $\mu
_{E}=\mbox{deg}E/\mbox{rk}E$. 
Two semistable bundles are called \emph{$S$-equivalent} if their
associated graded bundles 
\[
\mbox{Gr}(E)=\displaystyle\oplus _{i=1}^{m}E_{i}/E_{i-1} 
\]
are isomorphic, where 
\[
0=E_{0}\subset E_{1}\subset \cdots \subset E_{m}=E, 
\]
is the so-called Jordan-Holder filtration in which the successive quotients $%
E_{i}/E_{i-1}$ are stable of the same slope and are uniquely defined up to
permutation.
By the theorem of Narasimhan and Seshadri \cite{NS}, the space of $S$-equivalence 
classes of semistable bundles of rank $n$ and trivial determinant has the
structure of a projective algebraic variety, which we denote here by
${\cal M}_n(X)$.

In the case when $X=X_{\tau }=\Bbb{C}/(\Bbb{Z}\oplus \tau \Bbb{Z})$ is
an elliptic curve, and $E$ has trivial determinant, 
one can show that the Jordan-Holder quotients $L_{i}:=E_{i}/E_{i-1}$
are line bundles, necessarily of degree 0. Then $E$ is $S$-equivalent
to $L_1\oplus\cdots\oplus L_n$ and
\[
\mbox{det}E=L_{1}\otimes \cdots \otimes L_{n}=\mathcal{O}_{X_{\tau }}, 
\]
where $\mathcal{O}_{X_{\tau }}$ denotes the structure sheaf of $X_\tau$ 
(corresponding to the trivial line bundle on $X_\tau$).
Let $J(X_{\tau })\cong X_\tau$ be the Jacobian variety of $X_{\tau }$ and $M$ be
the kernel of the group homomorphism $t:J(X_{\tau })^{n}\to J(X_{\tau })$
given by the tensor product of line bundles. Consider the following 
natural maps
\[
\begin{array}{ccccc}
M\equiv \ker t & \to & \mathcal{M}_{n}(X_\tau) & \to & \mbox{Sym}^n (J(X_{\tau })) \\
(L_1,\dots , L_n) & \mapsto & E & \mapsto &  \{ L_{1},\dots , L_{n} \} 
\label{var.Abel.}
\end{array}
\]
where the last space is the symmetric product of the Jacobian.
One can prove that the second map is injective, and that the image of
the composition is a projective space. Therefore, one has
\begin{theorem}
{\emph{\cite{Tu}} For an elliptic curve $X_\tau$, the moduli space 
$\mathcal{M}_{n}=\mathcal{M}_{n}(X_\tau)$ is isomorphic to the complex projective space of dimension $n-1$, 
\[
\mathcal{M}_{n}\cong {\mathbb P}^{n-1}.
\]
Therefore, assuming $n\geq 2$, the Picard group $Pic(\mathcal{M}_{n})$ of
isomorphism classes of line bundles over }$\mathcal{M}_{n}${\ is isomorphic
to ${\mathbb Z}$, and letting $L_{\Theta }\cong{\cal O}(1)$ denote the ample generator, we
have 
\[
dimH^{0}(\mathcal{M}_{n},L_{\Theta }^{k})=\textstyle\binom{n+k-1}{k}.
\]
}
\end{theorem}

Although the analogous moduli spaces for higher genus curves do not admit
such a simple description, their Picard groups are all isomorphic to ${\ %
\mathbb Z}$ \cite{DN}. Therefore, in analogy with this case, a non-abelian theta
function (for genus 1) of level $k$ is defined to be a section $\theta $ of
the line bundle $L_{\Theta }^{k}$, 
\[
\theta \in H^{0}(\mathcal{M}_{n},L_{\Theta }^{k}). 
\]


\subsection{The Schottky Map}

\label{s3.2}

We now define the Schottky map for a complex linear subgroup $G$ of $GL(n,{%
\mathbb C})$ and a general compact Riemann surface $X$ of genus $g.$ Let $%
\mathcal{G}$ be the sheaf of germs of holomorphic functions from $X$ to $%
GL(n,{\mathbb C})$. The inclusion $G\hookrightarrow \mathcal{G}$ (where $G$
is identified with its constant sheaf on $X$), defines a map 
\begin{equation}
\mathcal{E}:H^{1}(X,G)\to H^{1}(X,\mathcal{G}),  \label{VT}
\end{equation}
that sends a flat $G$-bundle into the corresponding (isomorphism class of)
holomorphic vector bundle of rank $n$ over $X$ (necessarily of degree 0).

There is a well know bijection between the space of flat $G$-bundles over $X$
and the space of $G$-representations of the fundamental group of $X$, $\pi
_{1}(X)$, modulo overall conjugation $Hom(\pi _{1}(X),G)/G$; it is given
explicitely by: 
\begin{equation}
V:Hom(\pi _{1}(X),G)/G\to H^{1}(X,G),\qquad \rho \mapsto V_{\rho }:=\tilde{X}
\times _{\rho }{G},  \label{E.}
\end{equation}
where the notation means that $\pi _{1}(X)$ acts diagonally through $\rho $
on the trivial $G$-bundle over the universal cover $\tilde{X}$ of $X$.



Let us fix a canonical basis of $\pi _{1}(X)$: elements $%
a_{1},...,a_{g},b_{1},...,b_{g}$ that generate $\pi _{1}(X)$, subject to the
single relation $\prod_{i=1}^{g}a_{i}b_{i}a_{i}^{-1}b_{i}^{-1}=1$. Let $
F_{g} $ be a free group on $g$ generators $B_{1},...,B_{g}$, and $q:\pi
_{1}(X)\to F_{g}$ be the homomorphism given by $%
q(a_{i})=1,q(b_{i})=B_{i},i=1,...,g$. Then we can form the exact sequence of
groups: 
\begin{equation}
\pi _{1}(X)\stackrel{q}{\to }F_{g}\to 1,  \label{sequencia}
\end{equation}
which in turn defines the inclusion $i:Hom(F_{g},G)%
\hookrightarrow Hom(\pi _{1}(X),G)$.

\begin{definition}
The Schottky map is the composition $S=\mathcal{E}\circ V \circ i$, 
\[
G^{g}\cong Hom(F_{g},G)\stackrel{i}{\hookrightarrow }Hom(\pi _{1}(X),G)%
\stackrel{V}{\to }H^{1}(X,G)\stackrel{\mathcal{E}}{\to }H^{1}(X,%
\mathcal{G}).
\]

\end{definition}

Intuitively, this map sends a $g$-tuple of $n\times n$ invertible matrices 
\break
$(N_{1},...,N_{n})\in G^{g}\subset GL(n,{\mathbb C})^g$ to the flat rank $n$
holomorphic vector bundle determined by the holonomies 
$(1,...,1,N_{1},...,N_{g})$ around the loops 
\break
$(a_{1},...,a_{g},b_{1},...,b_{g}),$ respectively. To have a good description
of this map however, we need to substitute the space of all isomorphism
classes of vector bundles $H^{1}(X,\mathcal{G})$, by a nicer space such as
the moduli spaces of semistable bundles of section \ref{s3.1}. 
If $\rho\in Hom(F_g, U(n))$, it is known that the Schottky map 
is locally bi-holomorphic onto a neighbourhood in the moduli space of 
semistable bundles \cite{Fl}, however it is conjectured that its image is dense 
in ${\cal M}_n(X)$.  
In sections 
\ref{s3.3} and \ref{s3.4} we will consider two cases where this map can be
given explicitly.


\subsection{The rank one Schottky Map}

\label{s3.3}

In the case of line bundles the situation is simple, since the group $GL(1,%
\Bbb{C})=\Bbb{C}^{*}$ is an abelian group and the degree 0 line bundles in $%
H^{1}(X,\mathcal{O}^{*})$ form an abelian variety, the Picard variety of $X,$%
\[
Pic^{0}(X)=\frac{H^{1}(X,\mathcal{O})}{H^{1}(X,\Bbb{Z})}\cong \frac{H^{1}(X,%
\Bbb{C})}{H^{0}(X,K_X)\oplus H^{1}(X,\Bbb{Z})},
\]
where $K_X$ is the canonical bundle on $X$.
The last expression is obtained from the short exact sequence (see \cite{Gu}) 
\begin{equation}
0\rightarrow \Bbb{C}\rightarrow \Bbb{\mathcal{O}}\rightarrow K_X\rightarrow 0.  
\label{sec}
\end{equation}

Our choice of basis for $\pi _{1}(X)$ induces an isomorphism  $H^{1}(X,\Bbb{C}
)\cong Hom(\pi _{1}(X),\Bbb{C})\cong \Bbb{C}^{2g},$ and allows a very
explicit description of the Schottky map as follows 
\begin{equation}
\begin{array}{lccc}
S: & Hom(F_{g},\Bbb{C}^*)\cong (\Bbb{C}^*)^{g} & \Bbb{\rightarrow } & 
Pic^{0}(X)\cong \Bbb{C}^{2g}\diagup H^{0}(X,K_X)\oplus H^{1}(X,\Bbb{Z}) \\ 
& (e^{2\pi iz_{1}},...,e^{2\pi iz_{g}}) & \mapsto & 
[(0,...,0,z_{1},...,z_{g})].
\end{array}
\label{eq3.1}
\end{equation}

We can still be more concrete and at the same time
relate the Schottky map with the construction of \cite{FMN}, by using the
Jacobian variety $J(X)$ instead of the Picard variety of $X$.
By definition
\begin{equation}\nonumber
J(X)=H^{0}(X,K_X)^{*}/H_{1}(X,\Bbb{Z})\cong \Bbb{C}^{g}/\left( \Bbb{Z}
^{g}\oplus \Omega \Bbb{Z}^{g}\right), 
\end{equation}
where the last expression arises from considering 
$\{a_1,\dots,a_g,b_1,\dots,b_g\}$ as the basis of the lattice
$H_{1}(X,\Bbb{Z})$ and $\{a_1,\dots,a_g\}$ as the basis of 
the complex vector space $H^{0}(X,K_X)^{*}$, via the natural action of 
$H_{1}(X,\Bbb{Z})$ on $H^{0}(X,K_X)^{*}$.
As is well known, the Jacobian and
Picard varieties of $X$ are canonically isomorphic, as follows. Let 
\[
\Pi :H^{1}(X,\Bbb{C})\cong \Bbb{C}^{2g}\rightarrow H^{1}(X,\mathcal{O})\cong 
\Bbb{C}^{g} 
\]
be any linear map with kernel equal to $H^{0}(X,K).$ Simple computations
using (\ref{sec}), show that $\Pi $ is represented by a 
$g\times 2g$ matrix, also denoted $\Pi =\left[ \Pi
_{1}|\Pi _{2}\right] ,$ such that 
$ \Pi _{1}+\Pi_{2}\Omega =0.$
The isomorphism between the Picard and Jacobian
varieties of $X$ is then given by $\Pi _{2}:J(X)\rightarrow Pic^{0}(X),$ 
(see \cite{Gu}), hence 
\begin{equation}
\begin{array}{ccccc}
H^{1}(X,\Bbb{C}) & \stackrel{\Pi }{\rightarrow } & Pic^{0}(X) & \stackrel{
\Pi _{2}^{-1}}{\rightarrow } & J(X) \\ 
(w,z) & \mapsto & \Pi _{1}w+\Pi _{2}z & \mapsto & -\Omega w+z.
\end{array}
\label{eq3.3}
\end{equation}

Therefore, we obtain

\begin{theorem}
\label{t3.1} With our choices of basis of $H_{1}(X,\Bbb{Z})$ 
and $H^{0}(X,K_X)^{*}$, we have
\begin{equation}
\label{eq3.2}
(\Pi _{2}^{-1}\circ S)(e^{2\pi iz_{1}},\cdots ,e^{2\pi iz_{g}}) = z\,(\mathrm{mod} \Lambda ).
\end{equation}
\end{theorem}
\begin{proof}
This follows immediately from the maps (\ref{eq3.1}) and (\ref{eq3.3}).
\end{proof}

Because of this result, the composition $\Pi _{2}^{-1}\circ S:(\Bbb{C}%
^{*})^{g}\rightarrow J(X)$, is independent of the actual isomorphism between 
$J(X)$ and $Pic^{0}(X),$ and will henceforth be called the abelian Schottky
map, and denoted by $s$.
This map was used in \cite{FMN}, in order to identify, via pullback, classical theta
functions with holomorphic functions on $({\mathbb C}^{*})^{g}$.
We will use it later to relate the CST for the
elliptic curve $X_{\tau }$, in definition \ref{def41}, to non-abelian theta functions of genus one.

\subsection{The Schottky Map for genus one}

\label{s3.4}

Let us now consider the Schottky map for the case of an elliptic
curve $X_{\tau }=\Bbb{C}/(\Bbb{Z}\oplus \tau \Bbb{Z})$, and for semistable
bundles of rank $n$ with trivial determinant over $X_{\tau }$. This will be
a map 
\[
S:SL(n,{\mathbb C})\to \mathcal{M}_{n}\cong \Bbb{P}^{n-1}. 
\]
From geometric invariant theory (see \cite{GIT,N}), 
we know that under the adjoint action, $SL(n,{\mathbb C)}$ 
has a good quotient, which is a map $SL(n,{\mathbb C}
)\rightarrow T_{\Bbb{C}}/W$, where $T_{\Bbb{C}}$ is the maximal torus of $%
SL(n,{\mathbb C)}$ and $W$ the Weyl group; this map sends a matrix to the
unordered set of its eigenvalues. Therefore $SL(n,{\mathbb %
C)}$ satisfies a universal property, which in this case translates into the
statement that the Schottky map factors through $T_{\Bbb{C}}/W,$ as in the
following diagram. 
\[
\begin{array}{ccc}
SL(n,{\mathbb C}) & \stackrel{Q}{\to } & T_{\Bbb{C}}/W \\ 
& \searrow S & \downarrow f \\ 
&  & \mathcal{M}_{n}
\end{array}
\]
We are therefore reduced to describing $f.$ Let us fix the following Cartan
subalgebra of $sl(n,\Bbb{C})$
\[
\h=\{A\in SL(n,{\mathbb C):A}\text{ is diagonal}\}, 
\]
and using the coroot lattice $\check{\Lambda}_{R},$ define $M$ to be the
abelian variety 
\[
M=\check{\Lambda}_{R}\otimes X_{\tau }=\check{\Lambda}_{R}\otimes \left( 
\Bbb{C}/\Bbb{Z}\oplus \tau \Bbb{Z}\right) =\h/\left( \check{\Lambda}
_{R}\oplus \tau \check{\Lambda}_{R}\right) , 
\]
From the explicit form of $\h$ it is
not difficult to see that $M$ is isomorphic to the kernel of the map $%
t:J\left( X_{\tau }\right) ^{n}\rightarrow J\left( X_{\tau }\right) ,$
obtained by tensoring the entries, as in subsection \ref{s3.1}, with explicit
isomorphism as follows. 
\[
\begin{array}{ccc}
\h/\left( \check{\Lambda}_{R}\oplus \tau \check{\Lambda}_{R}\right) & 
\rightarrow & \ker t\subset J\left( X_{\tau }\right) ^{n} \\ 
\left[ (z_{1},...,z_{n})\right] & \mapsto & \left(
L_{z_{1}},...,L_{z_{n}}\right),
\end{array}
\]
where $z_{i}\in \Bbb{C},$ $,$ $z_{1}+\cdots +z_{n}=0,$ and $L_{z}$ denotes
the line bundle over $X_{\tau }=\Bbb{C}/(\Bbb{Z}\oplus \tau \Bbb{Z})$
corresponding to the divisor $[z]-[0].$ We will use these two different
representations interchangeably.

The Weyl group $W$ acts naturally on $M$, via the usual action on $\ck \La_R$, 
and, as shown in \cite{FM},\cite{FMW},\cite{L},\cite{M}, $\mathcal{M}_{n}$ is
isomorphic to the $l=(n-1)$-dimensional complex projective space $\Bbb{P}
^{l}(\Bbb{C})$ obtained as the quotient under this action,
\[
\mathcal{M}_{n}\cong M/W. 
\]
In our case, this quotient is given explicitely by the map \cite{Tu} 
\begin{equation}
\begin{array}{cccc}
\pi: & M & \rightarrow & \mathcal{M}_{n} \\ 
& (L_{z_{1}},...,L_{z_{n}}) & \mapsto & L_{z_{1}}\oplus \cdots
\oplus L_{z_{n}}.
\end{array}
\label{pi.2}
\end{equation}
Let now 
\[
s:{\mathbb C}^{*}\to J(X_{\tau })\cong X_{\tau },\qquad e^{2\pi iz}\mapsto
L_{z} 
\]
be the abelian Schottky map for genus one, where we identify, as usual the
elliptic curve with its Jacobian. It is easy to see that we have the
following commutative diagram, where the vertical arrows are inclusions 
\[
\begin{array}{ccccc}
\mathfrak{h} & 
\rightarrow & T_{\Bbb{C}} & \stackrel{s^{n}|_{T_{\Bbb{C}}}}{\rightarrow }
& M \\ 
\downarrow &  & \downarrow &  & \downarrow \\ 
\Bbb{C}^{n} & \rightarrow & \left( \Bbb{C}^{*}\right) ^{n} & \stackrel{s^{n}%
}{\rightarrow } &  J\left( X_{\tau }\right)^{n}.
\end{array}
\]
The map $f$ is then the Weyl invariant restriction of $s^{n}$ to $T_{\Bbb{C}%
} $. More precisely,

\begin{proposition}
\label{semlabel}
The following diagram is commutative 
\[
\begin{array}{cccclcl}
SL(n,{\mathbb C}) & \stackrel{Q}{\to } & T_{\Bbb{C}}/W & \stackrel{\tilde\pi }{
\leftarrow } & T_{\Bbb{C}} & \hookrightarrow  & \left( \Bbb{C}^{*}\right)
^{n} \\ 
& S\searrow  & \downarrow f &  & \downarrow s^{n}|_{T_{\Bbb{C}}} &  & 
\downarrow s^{n} \\ 
&  & \mathcal{M}_{n} & \stackrel{\pi}{\leftarrow } & M & 
\hookrightarrow  & J\left( X_{\tau }\right)^{n}
\end{array}
\]
\end{proposition}

\begin{proof}
The commutativity of the right square is clear, so let us consider
the middle one. By the construction of the good quotient, $f$ coincides with 
$S$ when evaluated on diagonal matrices 
\[
f\circ \tilde\pi (w_{1},...,w_{n})=f(\{w_{1},...,w_{n}\})=S(\mbox{diag}%
(w_{1},...,w_{n})),\qquad w\in T_{\Bbb{C}}. 
\]
Since the diagonal matrices correspond to vector bundles which are direct
sums of line bundles of degree 0, we have 
\begin{equation}
\label{nuncamaisacaba}
S(\mbox{diag}(w_{1},...,w_{n}))=s(w_{1})\oplus \cdots \oplus
s(w_{n})=L_{z_{1}}\oplus ...\oplus L_{z_{n}} 
\end{equation}
where $e^{2\pi iz_{j}}=w_{j},$ $j=1,...,n,$ and we identify once again $
X_{\tau }$ with $J(X_{\tau })$. Therefore, by (\ref{pi.2}), 
\[
f\circ \tilde \pi (w_{1},...,w_{n})=\pi(L_{z_{1}},...,L_{z_{n}})
=\pi\circ s^{n}(w_{1},...,w_{n}), 
\]
which proves the proposition.
\end{proof}



\section{Non-abelian Theta Functions in Genus One}
\la{s4}

\subsection{The Structure of Non-abelian Theta Functions in Genus One}
\la{s41}

Consider again our elliptic curve $X_\tau = \C / (\Z + \tau \Z)$ 
and the moduli space $\M_{n}(X_\tau)$ of
semistable holomorphic vector bundles of rank $n$
and trivial determinant over $X_\tau$. 
We are therefore restricting ourselves in this and next subsection  
to the case when $K=SU(n)$, with the extra assumption that $n\geq 3$. The case of $SU(2)$ is special and will 
be treated in subsection 
\ref{5.3}.  
The case
of general compact semi-simple group $K$ can be treated by analogous techniques. 
As above, let the Cartan subalgebra be 
${\h} = \{A \in sl(n,\C) \ : \ \hbox{\rm $A$ is diagonal} \} ,$
and consider the abelian variety
$$
M = \ck \La_R \otimes X_\tau = 
\h /  (\ck \La_R \oplus  \tau \ck \La_R). 
$$

Since our aim is to describe theta functions as
$Ad_{SL(n,\C)}$-invariant holomorphic functions
on $SL(n,\C)$
it will be usefull to recall their definition
as $W$-invariant sections of appropriate 
holomorphic line bundles on  $M$. 
Therefore, we apply the results of section \ref{s22} to the
abelian variety $M$.
All Weyl invariant antisymmetric integral forms $E$ on 
$\La =  \ck \La_R \oplus  \tau \ck \La_R$
are integral multiples of the form $E_1$
given by  \cite{L}
\begin{eqnarray}
\nonumber
& E_1 (\ck{\al}, \tau \ck{\beta}) & =  
<\ck{\al}, \ck{\beta}> = \, <\al,\beta>\\
& E_1 (\ck{\al}, \ck{\beta})      & =  
E_1 (\tau \ck{\al}, \tau \ck{\beta}) = 0, \quad \quad \forall \al, \beta \in \Lambda_R
\la{411}
\end{eqnarray}
where, as in section \ref{s21b}, $\ck \alpha, \ck \beta$
denote the coroots corresponding to $\al , \beta \in \La_R$ and 
recall that $<,>$ is the inner product on $\h$ for which the roots have squared length 2.

We are now interested in finding a classification of line bundles on $M$ which are Weyl invariant.
Recall that there exists a one to one correspondence between W-invariant antisymmetric integral bilinear forms
on $\Lambda$ which are compatible with the complex structure, and elements of the lattice of integral symmetric 
bilinear forms on $\ck\La_R$, denoted by $S^2\ck\La_R$ \cite{L}.
Consider the familiar exact sequence 
\begin{equation}
\la{casca2}
0 \rightarrow Pic^{0}(M) \rightarrow Pic(M) \stackrel{c_1}{\rightarrow} H^2(M,\Z).
\end{equation}
\begin{proposition} 
\la{casca}
[Looijenga \cite{L}]\footnote{The $Pic^0 (M)^W$ part of the sequence (3.2.1) of \cite{L} 
does not hold for $n\geq 3$ as we 
show in our proof. In subsection \ref{5.3} we show that it does hold for $n=2$.}
Let $n\geq 3$. For the abelian variety  $M = X_\tau \otimes \ck \La_R$ the sequence (\ref{casca2})
becomes
\begin{equation}
\la{casca1}
0 \rightarrow {\h}^* / (\La_W \oplus \tau \La_W) \rightarrow Pic(M) 
\rightarrow S^2 \ck \La_R \rightarrow 0,
\end{equation}
for which the Weyl invariant part is
\begin{equation}
\la{casca3}
0 \rightarrow 0 \rightarrow Pic(M)^W \rightarrow (S^2 \ck\La_R)^W \rightarrow 0.
\end{equation}
\end{proposition}
\begin{proof}
We only show that $Pic^{0}(M)^W=0$ since this is the only difference with respect to \cite{L}.
The automorphy factors of the line bundle corresponding to 
$x\in \h^*$, with $x=x_1 +\tau x_2$ and $x_1,x_2 \in \La_R \otimes \R = \h_\R^*$ 
are given by 
\begin{equation}
\la{casca4}
e(\ck \al + \tau \ck \beta) = e^{2\pi i (x_1(\ck \al) + x_2(\ck \beta))}, \,\,\,\,\, \ck \al, \ck 
\beta \in \ck\La_R.
\end{equation}
So, $x$ leads to a Weyl invariant point in $Pic^{0}(M)$ if and only if for all the elementary Weyl 
reflections $w_j$ one has 
\begin{equation}
\la{casca5}
x-w_j(x) = <\al_j,x>\al_j \in \La_W\oplus\tau\La_W.
\end{equation} 
{}From the Cartan matrix for the algebras $A_l$ for $l\geq 2$ we conclude that 
this implies that $<\al_j,x> \in \Z \oplus \tau\Z$, for all simple roots $\al_j$, 
and therefore $x \in \La_W\oplus\tau\La_W$ and leads to the trivial line bundle $0 \in Pic^{0}(M)$.
\end{proof}

Since all symmetric Weyl invariant integral bilinear forms on $\ck \La_R$ are 
integral multiples of $<,>$, we see from 
(\ref{casca3}) that $Pic(M)^W$ is infinite cyclic $i.e.$ there is a Weyl invariant line bundle 
$L_1 \rightarrow M$ (this is the same line bundle $L_1$ of section \ref{s22}, since $c_1(L_1)$ 
corresponds to $E_1$) such that 
\begin{equation}
\la{casca6}
Pic(M)^W = \{ L_1^k, k\in \Z\}.
\end{equation}

Recall that we have the projection $\pi: M \rightarrow \M_n(X_\tau)$, and let 
$L_\Theta$ be the theta bundle 
over $\M_n(X_\tau)$, as in section \ref{s3}.
{}From theorem (3.4) of \cite{L}, we conclude that $\pi^* L_\Theta^k \cong L_1^k$.
More precisely, we have
\begin{lemma}
\la{casca88}
Let $\pi :M\rightarrow M/W\cong \P^{n-1}$ be the canonical projection and 
${\cal L}=\pi^{*}L_\Theta.$ Then ${\cal L}=L_{1},$ for $n\geq 3$.
\end{lemma}
\begin{proof}
Since ${\cal L}$ is a $W$-invariant line bundle on $M,$ its polarization 
$E=-{\rm Im}H$ is a multiple of $E_{1},$ which means that ${\cal L}=L_{1}^{p},$ for
some $p\in \Z,$ $p\geq 1.$ 
It is shown in \cite{L} that 
\[
H^{0}(M,L_1)^{W}\cong \C^{n},
\]
and that $\dim H^0(L_1^m)^W >n$ for $m>1$. Since $L_\Theta \cong {\cal O}(1)$ we have 
$n=h^{0}(\P^{n-1},L_\Theta)=\dim H^0({\cal L})^W$ which implies ${\cal L}=L_{1}$.
\end{proof}

In the notation of section \ref{s22} we have that 
$\La_1 = \ck \La_R$, $\La_2 = \tau \ck \La_R$ and  
$\wh \La = \ck \La_W \oplus \tau \ck \La_W$. The form $F$ on $V\times V_2$ is given by 
\begin{equation}
F (.,.)= - \tau^{-1} <.,.>. 
\label{412} 
\end{equation}
The space of theta functions 
$H^0(M, L_1^k)$ is isomorphic to a subspace of the space 
of holomorphic functions 
on $\h /\ck\La_R \cong (\C^*)^{l} \cong T_\C$ with a basis 
$\{{\theta_{\gamma,k}}\}_{\gamma\in \La_W/k\La_R}$ 
given by (see  (\ref{2.2.9})), 
\begin{equation}
\label{413}
\theta_{\gamma , k}(v) =  \sum_{\al\in \La_R}    
e^{\pi i k \tau <\al +\frac{\gamma}{k}, \al +\frac{\gamma}{k}> 
+ 2\pi i k (\al +\frac{\gamma}{k})(v)}, \quad v\in\h.
\end{equation}

As shown in section  \ref{s22},  these functions are the image 
under an abelian CST of certain distributions in $U(1)^l$. 
In order to apply theorem \ref{t23.1} to (\ref{413}), we show 
in the appendix that the basis of simple coroots $\ck \al_i$, 
for $i=1,...,l$ cannot be completed to a canonical basis of 
$\La_1\oplus \La_2 = \ck\La_R \oplus \tau \ck\La_R$ so that we are in the 
situation described in (\ref{331})-(\ref{336}). The role of 
$\gamma'_{i+l}$ is being played by $\tau \ck \lambda_i$, where $\ck \lambda_i$ are 
the fundamental coweights. The 
matrix $\tilde \Omega$ is then given by $\tilde \Omega = \tau C^{-1}$ 
where $C^{-1}$ is the inverse Cartan matrix and $R=C$ (see (\ref{331})).  
The automorphy factors read   
\begin{equation}
\la{413.5}
e(\tau\ck \al,v) = e^{-2\pi i k \al (v) - \pi i k \tau 
<\al , \al>}
\end{equation}
in agreement with formula (\ref{336}).

The pull-back under $\pi \circ s^n|_{T_\C}$ of  $H^0 ({\cal M}_{n}(X_\tau), L_\Theta^k)$ 
to ${\h}/\ck\La_R$, 
corresponds the space $\H_{k,\tau}^+$ of  Weyl invariant linear combinations of elements
of the form (\ref{413}).
The spaces of level $k$ non-abelian theta functions $\H_{k,\tau}^+$ are
the fibers of a vector bundle over the Teichm\"uller space 
of genus one curves (called the bundle of genus one, level $k$, $SU(n)$
conformal blocks in conformal field theory)
\ba
\la{4121}
\H_k \ & \longrightarrow & \ {\cal T}_1 = \left\{\tau \in \C \ : \ {\rm Im} 
\tau = \tau_2 >0 \right\} = \HH_1 \\
\H_{k_{|_\tau}} &=& \H_{k,\tau}^+ \ .  \nonumber
\ea
{}From (\ref{413}) we see that, for every $\tau$,  a basis of  
$\H_{k,\tau}^+$ is given 
by  Weyl invariant theta functions of the form 
\begin{equation}
\label{414}
\theta_{\gamma ,k}^+ = \sum_{w\in W} \theta_{w (\gamma ), k},   
\qquad  \gamma \in  \La_W/ (W \triangleright k\La_R),
\end{equation}
where $W \triangleright k\La_R$ denotes de semi-direct product of $W$ and $k\La_R$.

Taking into account the $\tau$ dependence, $\{\theta_{\gamma ,k}^+\}$ defines
a global moving frame of sections of $\H_k \rightarrow {\cal T}_1$
and therefore fixes a trivialization of the bundle of conformal
blocks. A different trivialization is obtained as follows.
 Let $\H_{k,\tau}^-$ be the space of Weyl anti-invariant theta functions of 
level $k$ with basis given by
\begin{equation}
\label{415}
\theta_{\gamma ,k}^- = \sum_{w\in W} \epsilon (w)\,\,\,\theta_{w (\gamma ), k},  \qquad \gamma \in 
 \La_W/ (W \triangleright  k\La_R)   \ .
\end{equation}  
Notice that $\theta_{\gamma ,k}^- = 0$  if $\gamma$  is singular, i.e. if $<\gamma,\alpha_i>= 0$  
for some simple root $\alpha_i$. A non-singular dominant weight 
$\gamma \in \La_W^+$ can always be written in the form 
$\gamma = \gamma' + \rho$ with $\gamma'$ being a, possibly singular, dominant weight.

Recall the following 
\begin{theorem} 
\label{tl}
[Looijenga \cite{L}]
\begin{itemize}
\item[a)] The space of Weyl anti-invariant theta functions of level $n$,  $\H_{n,\tau}^-$ is 
one-dimensional and $\H_{n,\tau}^- =   <\theta_{\rho ,n}^->_{\C}$.
\item[b)]
The map 
\begin{eqnarray}
\label{416}
\nonumber
\H_{k,\tau}^+ & \rightarrow & \H_{k+n,\tau}^- \\
\theta^+ & \mapsto & \theta^- = \theta_{\rho ,n}^- \ \theta^+ ,
\end{eqnarray}
is an isomorphism between the space of level $k$ Weyl invariant theta functions and the space 
$\H_{k+n,\tau}^-$ of Weyl anti-invariant theta functions of level $k+n$. 
\end{itemize}
\end{theorem}

Let $D_k$ be the set of integrable representations of level $k$ of the Kac-Moody algebra $sl(n,\C )_k$,  
$D_k = \{ \lambda\in \La_W^+ | \,\,\,< \lambda, {\hat \alpha}> \,\,\,\,\leq k\} 
\cong \La_W / (W \triangleright k \La_R)$, where $\hat \al = \al_1+\cdots +\al_{n-1}$ is the highest root for 
$sl(n,\C )$. 
The previous theorem leads to a  basis of $\H_{k,\tau}^+$, different from (\ref{414}), given by
\begin{equation}
{\hat \theta}^+_{\gamma,k} = \frac{\theta^-_{{\gamma}+\rho, k+n}}
{\theta^-_{\rho,n}}, \qquad {\gamma}\in D_k. 
\label{418}
\end{equation}
As 
we will see in the next subsection
from the point of view of the heat equation 
the trivialization of the bundle of 
conformal blocks corresponding to (\ref{418})
is more convenient than the one corresponding to 
(\ref{414}).


\subsection{CST and Non-abelian Theta Functions in Genus One}
\la{s42}

We now continue to follow the strategy indicated in the introduction.
As we mentioned before, by extending the $SU(n)$-CST $C_t^{\tau}$ of definition \ref{def41}
to distributions we of course lose unitarity.
From propositions \ref{p41} and \ref{p42} it follows that
the image under $C_t^{\tau}$ of a distribution on
$SU(n)$ with infinite $L^2$ norm is a holomorphic
function on $SL(n,\C)$ with infinite norm
with respect to the heat kernel measure $d\nu_{t\tau_2}$.

Consider again the projection 
$$Q: SL(n,\C) \rightarrow T_\C/W \cong {\h}/(W \triangleright \ck \La_R) \cong SL(n,\C)/Ad_{SL(n,\C)}$$
of section \ref{s3.4}.
In order to recover unitarity we will, for every $\tau$,
restrict integration on $SL(n,\C)$ to the
region $Q^{-1}([\h_0])$ with $\h_0 \subset \h$
a fundamental domain with respect to the group
$W \triangleright (\ck \La_R \oplus \tau \ck \La_R)$ and $[\h_0]= W\h_0 + \ck \La_R$. 
Of course this
will be meaningful only for those holomorphic 
funtions on $SL(n,\C)$ for which the integral
will not depend on the choice of $\h_0$. As we will see, this simple condition gives a precise 
analytic characterization of the pull-back of non-abelian theta functions to $SL(n,\C)$ with 
respect to the Schottky 
map $S$ (theorems \ref{levico6} and \ref{levico8}). The relation between the heat kernel
measures $d\nu_{t\tau_2}$ on $SL(n,\C)$ and $d\nu_{t\tau_2}^{\rm ab}$ on $T_\C$ is given by
\begin{proposition}
\label{p5}
The push-forward of the measure $d\nu_{t\tau_2}$ on $SL(n,\C)$ with 
respect to the projection $Q$ is given by 
\begin{equation}
\label{501b}
Q_{*} d\nu_{t\tau_2} = e^{-2t\pi \tau_2 ||\rho||^2} {|\sigma|^2} 
d\nu_{t\tau_2}^{\rm ab}
\end{equation}
where on the right-hand side we denote the restriction of $d\nu_{t\tau_2}^{\rm ab}$ to a fundamental domain of 
$W \triangleright \ck \La_R$ in $\h$, by the same symbol.
\end{proposition}

\begin{proof}
This follows from the fact that $\varphi_\C$ in (\ref{b11}) is an 
isometry.
\end{proof}

Consider then the following problems. 
We will use the same notation for $(W\triangleright \ck \Lambda_R)$-invariant functions 
on $\h$, the corresponding $Ad$-invariant functions on $SL(n,\C)$ and also 
their restrictions to $T_\C$. 

\begin{problem}
\label{problem1} 
Find $t>0$ such that 
there exist non-trivial subsets of 

$\left(C^{\infty}(SU(n)\right)^\prime)^{Ad_{SU(n)}}$ 
defined by  
\ba
\nonumber
{\cal F}  = && \{  \psi\in \left(C^{\infty}(SU(n)\right)^\prime)^{Ad_{SU(n)}}:  
 |C_t^{\tau}(\psi)|^2  |\sigma|^2   \nu_{t\tau_2}^{\rm ab}, \\
\label{502}
&& {\rm is}  \,\,\,\, \tau \ck\La_R \,\,\, {\rm invariant\, as\, a\, function\, on}\, \h \}.
\ea
\end{problem}

\begin{remark}
\label{remark1}
Since the functions $|C_t^{\tau}(\psi)|^2  |\sigma|^2   \nu_{t\tau_2}^{\rm ab}$ in (\ref{502}) are automatically 
$W \triangleright \ck\La_R$ invariant, being also $\tau \ck \La_R$ invariant means that they are  
the pull-back of functions on 
${\cal M}_n = \h / (W \triangleright (\ck \La_R \oplus \tau \ck \La_R))$. 
\end{remark}

For $\Psi$ and $\Psi'$ $Ad$-invariant holomorphic functions on $SL(n,\C)$,
and for a choice of a fundamental domain ${\h}_0 \subset \h$ of the group 
$W \triangleright  (\ck \La_R \oplus \tau \ck \La_R)$, define the modified Hall inner product as
\begin{eqnarray}
\label{503z}
\label{831}
\langle\langle\, \Psi, \Psi' \,\rangle\rangle := \int_{Q^{-1}([{\h}_0])} 
\overline{ \Psi}  \Psi' d\nu_{t\tau_2} = \int_{{\h}_0} \overline{ \Psi} \Psi' 
e^{-2t\pi\tau_2 ||\rho||^2} |\sigma|^2 d\nu_{t\tau_2}^{\rm ab},
\end{eqnarray}
where $[{\h}_0] = W{\h}_0 + \ck \Lambda_R$.

\begin{problem}
\label{problem2}
Find $t>0$ such that the modified Hall inner product 
\begin{eqnarray}
\label{503}
\langle\langle\, C_t^\tau (\psi), C_t^\tau (\psi')\,\rangle\rangle 
\end{eqnarray}
is independent of $\tau$ for all $\psi,\psi'$ belonging to ${\cal F}$. 
\end{problem}

We will see that, remarkably, the spaces $\cal F$ providing solutions to problem \ref{problem1} are
 independent of $\tau$ and   
lead to functions $C_t^\tau(\psi) \sigma$ 
which solve problem  \ref{problem2} and are automatically the pull-backs to $\h/\ck \La_R$ of holomorphic sections of line bundles over 
$\h/(\ck \La_R \oplus \tau \ck\La_R)$.

We know from lemma 4.2 of \cite{FMN} and section \ref{s22} above that if $t=1/k$ with 
$k\in\mathbb N$ then $\nu_{t\tau_2}^{\rm ab}$ defines an 
hermitean structure on $L_1^k$. 
We then have the following result,

\begin{lemma}
\label{levico1}
Problem \ref{problem1} has no solution if $t \neq 1/k'$ with $k' \in {\mathbb N}$.
\end{lemma}

\begin{proof}
The key observation is that $\nu_{t\tau_2}^{\rm ab}$ satisfies the quasi-periodicity conditions
\ba
\nonumber
\nu_{t\tau_2}^{\rm ab}(v +  \ck \al) & = &  \nu_{t\tau_2}^{\rm ab}(v),\\
\la{507}
\nu_{t\tau_2}^{\rm ab}(v + \tau \ck \al) & = &
|e^{2\pi i  \frac{1}{t} \al(v)}|^2|e^{\pi i  \frac{1}{t} \tau <\al,\al>}|^2 \nu_{t\tau_2}^{\rm ab}(v), 
\,\, {\rm for \, all}\,\ck \al \in \ck \La_R, 
\ea
which follow from (\ref{iproclaim}). On the other hand, from corollary \ref{co41} $C_t^\tau (\psi)\sigma$ is a holomorphic 
function on $\h \times {\cal T}_1$ and verifies 
\begin{equation}
\la{508}
\left(C_t^\tau (\psi) \sigma \right)(v + \ck \al_j)  =   \left(C_t^\tau (\psi) \sigma \right)(v),
\end{equation}
for all $\psi\in \left(C^{\infty}(SU(n)\right)^\prime)^{Ad_{SU(n)}}$. 

If the function $|C_t^{\tau}(\psi)|^2  |\sigma|^2   \nu_{t\tau_2}^{\rm ab}$ is 
$(\ck\La_R \oplus \tau \ck\La_R)$-invariant then from (\ref{507}) and the holomorphicity in 
$\h \times {\cal T}_1$ it follows that $C_t^\tau (\psi)\sigma$ must also satisfy the following 
quasi-periodicity conditions
\begin{equation}
\label{levico2}
\left(C_t^\tau (\psi) \sigma \right)(v + \tau \ck \al)  =  e^{- 2\pi i <d,\al>}
e^{-2\pi i \frac{1}{t} \al(v)}e^{-\pi i  \frac{1}{t} \tau <\al,\al>} 
\left(C_t^\tau (\psi)\sigma \right)(v), 
\end{equation}
where $d\in \h_\R^*$. 
Non-zero holomorphic functions satisfying (\ref{508}) and (\ref{levico2}) do not exist if
$1/t \notin {\mathbb N}$. This results from the fact that if $1/t \notin {\mathbb N}$ 
then the automorphy factors in (\ref{levico2}) are not invariant under $v \mapsto v+\ck\beta$, $\ck 
\beta \in \ck\Lambda_R$, 
which makes it impossible to solve (\ref{508}) and (\ref{levico2}).   
\end{proof}

Note that functions satisfying (\ref{508}) and (\ref{levico2}) with $t=1/k'$ are level $k'$ theta 
functions with automorphy factors (compare with (\ref{336}) and (\ref{413.5}))
\begin{equation}
\label{509}
e(\ck \al + \tau\ck \al',v) = e^{-2\pi i k' \al'(v) - \pi i k' \tau < \al',  \al'> - 2\pi i <d, \al'>}.  
\end{equation}

Let us now find the nontrivial set of distributions in 
$(C^\infty(SU(n))')^{Ad_{SU(n)}}$ solving problem \ref{problem1}. Since $\sigma$ is $W$-anti-invariant and 
$C^\tau_{1/k'}(\psi)$ is $W$-invariant, we conclude that the functions $C^\tau_{1/k'}(\psi)\sigma$ 
satisfying (\ref{508}) and (\ref{levico2}) are $W$-anti-invariant. 
As in the proof of proposition \ref{casca}, this implies that $d$ in (\ref{509}) belongs to $\La_W$.  
Therefore, we can take $d=0$. 

We then have,
\begin{theorem}
\label{levico6}
Problem \ref{problem1} has a solution, denoted by ${\cal F}_k$ 
if and only if $t=\frac{1}{k+n}$ with $k\in {\mathbb N}\cup \{0\}$.
The spaces ${\cal F}_k$ are $\tau$-independent, have dimensions given by the Verlinde numbers 
\begin{equation}
\la{503c}
{\rm dim\ }{\cal H}_{k+n,\tau}^{-} = \textstyle\binom{n+k-1}{k}
\end{equation}
 and have basis formed by 
\begin{equation}
\la{510}
\psi_{\gamma,k}(v) = \frac{1}{\sigma} \sum_{w\in W} \epsilon (w) \theta^0_{\gamma +\rho,k+n} (w(v)),\quad v\in \h_\R,
\end{equation}
where $\gamma =0$ if $k=0$ and $\gamma \in \La_W / (W \triangleright k\La_R)$ if $k>0$ and 
\begin{equation}
\la{511}
\theta^0_{\gamma+\rho,k+n} (v) = \sum_{\al \in \La_R} e^{2\pi i (\gamma + \rho + (k+n)\al)(v)}. 
\end{equation}
\end{theorem}

\begin{proof}
{}From lemma \ref{levico1} we can consider $1/t \in {\mathbb N}$. On the other hand, 
from theorem \ref{tl} it follows that non-zero $W$-anti-invariant theta functions exist only 
for level $k' \geq n$. 
So if $0\neq \psi \in (C^\infty(SU(n))')^{Ad_{SU(n)}}$ 
satisfies the condition (\ref{502})
 then $t=1/(k+n)$, $k\geq 0$ and  $\varphi_\C \circ C_{1/(k+n)}^{\tau}(\psi)$ 
(see (\ref{car2}))
is a Weyl
anti-invariant theta function of level $k+n$.  
Conversely, from theorems \ref{pcar} and \ref{t23.1} 
it follows that every Weyl anti-invariant theta function $\theta \in {\cal H}_{k+n,\tau}^{-}$ is the image under 
$\varphi_\C \circ C_{1/(k+n)}^{\tau}$ of a unique $Ad_{SU(n)}$-invariant distribution from ${\cal F}_k$.
It is easy to check that the inverse images 
under $\varphi_\C \circ C_{1/(k+n)}^{\tau}$ of the  theta functions $\theta^{-}_{\gamma,k}$ in (\ref{415})  
are given by the distributions $\sqrt{|W|}\psi_{\gamma,k}$ in the theorem. 
The $\tau$-independence of ${\cal F}_k$ follows 
from 
(\ref{510}).
\end{proof}

In terms of characters of irreducible representations of $SU(n)$, the distributions $\psi_{\gamma, k}, 
\gamma \in D_k = \{ \lambda\in \La_W^+ | \,\,\,< \lambda, {\hat \alpha}> \,\,\,\,\leq k\} 
\cong \La_W / (W \triangleright k \La_R)$, where $\hat \al$ is the highest root for 
$sl(n,\C )$,  
have the form 
\begin{equation}
\la{tasol}
\psi_{\gamma, k} = \sum_{\stackrel{\lambda \in \La_W^+}{\lambda+\rho 
\in [\gamma+\rho]}} \epsilon_\lambda \chi_\lambda ,
\end{equation}
where $[\gamma+\rho]$ is the $(W\triangleright (k+n) \Lambda_R)$-orbit of $(\gamma+\rho)$ and 
$\epsilon_\lambda = 
\epsilon(w)$ for the unique $w\in W$ such that $\lambda +\rho = w (\gamma +\rho)\,\, 
{\rm mod}\,\, (k+n)\La_R$.
This follows from the Weyl character formula (\ref{b7}) and the fact that the group $W \triangleright (k+n)\La_R$ 
acts freely on the set of dilated
Weyl alcoves \cite{Ber,Fr,PS}.

We now consider the images of the spaces ${\cal F}_k$ under the CST and show that the solution of problem \ref{problem1} 
leads to the solution of problem \ref{problem2}. 
Let $\wt \H_{k,\tau}$ be
the image of ${\cal F}_k$ under 
$C_{t=1/(k+n)}^{\tau}$,
\begin{equation}
\la{22-1}
\wt \H_{k,\tau} = C_{t=1/(k+n)}^{\tau}\left({\cal F}_k\right) . 
\end{equation}
On $\wt \H_{k,\tau}$ we have the hermitean inner product   $\langle\langle \cdot , \cdot \rangle\rangle$ 
(\ref{831}) induced by the Hall CST inner product.
The hermitean structure on 
$\H_k$ in (\ref{4121}), which is of interest in Conformal Field Theory, is \cite{AdPW}  
\begin{equation}
\label{506}
<\theta^+,\theta^{+ '}> = \frac{1}{|W|} \int_M \overline{ \theta^+} \theta^{+ '} 
|\theta^-_{\rho,n}|^2 d\nu_{\tau_2 / (k+n)}^{\rm ab}, 
\end{equation}
for $\theta^+, \theta^{+ '}\in \H_{k,\tau}^+$. We then have,
\begin{theorem} 
\label{levico8}
The family $\{\wt \H_{k,\tau}\}_{\tau \in {\cal T}_1}$ in (\ref{22-1}) forms a vector 
bundle over ${\cal T}_1$ isomorphic to the bundle of conformal blocks
$\H_k \rightarrow {\cal T}_1$ (\ref{4121}).

The hermitean structure defined by (\ref{831}), with $t=1/(k+n)$, does not depend on $\h_0$ and the  
map $\Phi_k : \wt \H_k \rightarrow \H_k$ in (\ref{4220}) is a unitary 
 isomorphism of vector bundles.
\end{theorem}

\begin{proof}
We have shown in theorem \ref{levico6} that the restriction of $\varphi_\C \circ C_{1/(k+n)}^{\tau}$
to ${\cal F}_k \subset (C^\infty(SU(n))')^{Ad_{SU(n)}}$ is an isomorphism to
$\H^-_{k+n,\tau}$. From theorem \ref{tl} it then follows that the bundle of conformal blocks
${\cal H}_k \rightarrow {\cal T}_1$ in (\ref{4121}) is isomorphic to the bundle  
\begin{equation}
\la{pp1}
\wt \H_{k} = \left\{ \wt \H_{k,{\tau}} \right\}_{\tau \in {\cal T}_1}  \rightarrow {\cal T}_1  
\end{equation}
with simple isomorphism $\Phi_k$ given by (\ref{4220}). 

The identities (\ref{507}), (\ref{508}) and (\ref{levico2}) imply that the hermitean structure 
(\ref{831}) does not depend on $\h_0$.
From theorems \ref{t23.1} and \ref{tl} it follows that (\ref{506}) defines a 
shifted hermitean structure on 
$\H_{k}$ for which the frame 
$$\{\wh \theta_{\gamma,k}^+ = \frac{\theta^-_{\gamma + \rho,k+n}}{\theta^-_{\rho, n}} \}$$ 
is orthonormal. (Notice that the same is not true for the ``unshifted'' frame $\{\theta_{\gamma,k}^+ \}$.) 
We see from formula 
(\ref{831}) that with these hermitean structures on $\wt \H_{k}$ and $\H_{k}$ the isomorphism $\Phi_k$ in 
(\ref{4220}) is unitary.
\end{proof}

Finally, we recover the unitarity of the CST with 

\begin{theorem}
\label{levico9} 
With respect to the $\tau$-independent inner product $\left(\cdot , \cdot \right)$ on ${\cal F}_k$ for 
which the basis $\{\psi_{\gamma,k}\}$ in 
(\ref{510}) and (\ref{tasol}) 
is orthonormal, 
the CST  
$$
C_{t=1/(k+n)}^\tau : \quad \left( {\cal F}_k, \left( \cdot, \cdot \right) \right) \rightarrow 
\left( \wt \H_{k,\tau}, \langle\langle\,\cdot ,\cdot \,\rangle\rangle \right)
$$ 
is a unitary isomorphism 
$\forall\tau\in {\cal T}_1$ and $k\in\N_0$.
\end{theorem}

\begin{proof}
{}From theorem \ref{pcar}  we see that 
\begin{equation}
\label{fim1}
C^{\tau}_{1/(k+n)} (\psi_{\gamma, k}) = \frac{e^{-\frac{i\pi \tau}{k+n}||\rho||^2}}{\sigma} \theta^-_{\gamma +\rho, k+n}.
\end{equation}
Unitarity then follows from (\ref{831})  and theorem \ref{t23.1}.
\end{proof}

In terms of characters of $SU(n)$, from (\ref{426}), (\ref{tasol}) and (\ref{fim1}) we obtain
\begin{equation}
\label{queijodeovelha}
C^{\tau}_{1/(k+n)} (\psi_{\gamma, k}) = \sum_{\stackrel{\lambda \in \La_W^+}{\lambda +\rho \in [\gamma+\rho]}} \epsilon_\lambda 
e^{\frac{i\pi \tau}{k+n} c_\lambda} \chi_\lambda .
\end{equation}


\subsection{The Case of $SU(2)$}
\label{5.3}

Somewhat surprisingly, the case of $K = SU(2)$ is special for several reasons. 
To begin with, the Weyl invariant antisymmetric integral bilinear forms on 
$$
\La = \ck \La_R \oplus \tau \ck \La_R = \Z \ck\al_1 \oplus \tau \Z \ck\al_1
$$ 
are integral multiples of $E_1$, where 
$E_1(\ck \al_1, \tau \ck\al_1)=\frac{1}{2}<\al_1,\al_1> = 1$ which is different from 
(\ref{411}) by a factor of $1/2$.  The form $F$ is now given by $F (.,.)= -\frac{1}{2}\tau^{-1} <.,.>$. 
As in (\ref{2.2.6}) this defines the line bundles $L_1^ k$ on the abelian variety 
$M = X_\tau \otimes \ck \La_R = \C \ck\al_1/ (\Z\ck\al_1 \oplus \tau \Z \ck \al_1)\cong X_\tau, 
 \C \ck\al_1 = \h$. 

The space of theta functions 
$H^0(M, L_1^k)$ is isomorphic to a subspace of the space 
of holomorphic functions 
on $\h /\ck\La_R \cong \C^* \cong T_\C$ with a basis 
$\{{\theta_{\gamma,k}}\}_{\gamma\in \La_W/(k/2)\La_R}$ 
given by (see  (\ref{2.2.9})), 
\begin{equation}
\label{531}
\theta_{\gamma , k}(v) =  \sum_{\al\in \La_R}    
e^{\pi i k \tau \frac{1}{2}<\frac{2}{k}\gamma +\al, \frac{2}{k}\gamma +\al> 
+ \pi i k (\frac{2}{k}\gamma +\al)(v)},
\end{equation}
or in more explicit classical notation, $\gamma = m \lambda_1 = (m/2) \al_1$, $v = z \ck\al_1$ and $\al = p\ck \al_1$
\begin{equation}
\la{maio1}
\theta_{\gamma,k}(v) = \theta_{m,k}(z) = \sum_{p\in \Z} e^{\pi i \frac{\tau}{k}(m+kp)^2 + 2\pi i (m+kp)z},
\end{equation}
with $0\leq m < k, \ m \in \N$.
In this case we have $(S^2\ck\La_R)^W = S^2\ck\La_R$ and 
\begin{proposition}
\la{casca444}
[Looijenga \cite{L}]
The Weyl invariant part of the sequence (\ref{casca1}) is
\begin{equation}
\label{casca445}
0 \rightarrow \left[1/2\left(\La_W \oplus \tau\La_W\right)\right]/\left(\La_W \oplus \tau\La_W\right) 
\rightarrow (Pic(M))^W 
\rightarrow S^2\ck\La_R \rightarrow 0
\end{equation} 
\end{proposition}
\begin{proof}
As in proposition \ref{casca} it suffices to determine $Pic^0(M)^W$.
The condition (\ref{casca5}) has a solution if and only if 
$x= (m_1 + \tau m_2) \alpha_1 /4$, $m_1,m_2 \in \Z$. Therefore, 
$$
\la{casca446}
\nonumber
(Pic^0 (M))^W \cong \left[1/2\left(\La_W \oplus \tau \La_W\right)\right]/ 
\left( \La_W \oplus \tau \La_W\right) \cong 
\Z_2 \oplus \Z_2.
$$\end{proof}

The second main difference is that, in contrast with lemma \ref{casca88}, 
not all line bundles in $Pic(M)^W$ are the pull-backs of line bundles on $\M_2(X_\tau)$. 
Denote line bundles in $(Pic^0(M))^W$ by $N_\eta$, $\eta \in \left[1/2\left( \La_W \oplus \tau \La_W\right)\right]
/ \left( \La_W \oplus \tau \La_W\right) $ and 
let $L_1$ be the line bundle on $M$ with automorphy factors defined by the form $F$ above,
\begin{equation}
\la{casca 447}
e(\ck\alpha + \tau \ck\beta, v) = e^{-\pi i \beta(v) - \frac{1}{2}\pi i \tau <\beta,\beta>}.
\end{equation} 
We see from (\ref{casca445}) that 
\begin{equation}
\la{casca448}
(Pic(M))^W = \left\{ N_\eta \otimes L_1^k, k\in \Z, \eta \in \left[1/2\left(\La_W \oplus \tau \La_W\right)\right]/ 
\left( \La_W \oplus \tau \La_W\right)  
\right\}.
\end{equation}
\begin{lemma}
\la{casca449}
Let $\pi: M \rightarrow M/W \cong {\cal M}_2(X_\tau)\cong \P^1$ be the canonical projection and 
${\cal L} = \pi^* L_\Theta$. 
Then ${\cal L} \cong L_1^2$.
\end{lemma}
\begin{proof}
In this case, $M\cong X_\tau$ and  the map $\pi:M \rightarrow {\cal M}_2$ becomes the 
two-fold ramified covering of 
$\P^1$, $X_\tau \rightarrow X_\tau/\Z_2 \cong \P^1$. The pull-back to $X_\tau$ of the 
the $k$th power of the theta bundle, ${\cal O}(k) \rightarrow \M_2$, will be a line bundle 
of degree $2k$. 
On the other hand, $\pi^{-1}([\pi((\frac{1}{2}+  \tau\frac{1}{2} )\ck\alpha_1)]) = 
2 [(\frac{1}{2}+  \tau\frac{1}{2} )\ck\alpha_1]$ as an element in $Div(X_\tau)$ and therefore 
$\pi^* {\cal O}(1)\cong L_1^2$, since the zero of the Riemann theta function 
($k=1, m=0$ in (\ref{maio1})) is $(1/2 +1/2 \tau)\ck\al_1$.
\end{proof}

As a consequence of the above results, we obtain as the third significative 
difference with the case $n\geq 3$, the fact that 
the CST and problems \ref{problem1} and \ref{problem2} lead not only to 
non-abelian theta functions on ${\cal M}_2(X_\tau)$, but also to Weyl invariant 
theta functions on $M$ which do not descend to sections of bundles on the moduli space.
We have in place of lemma \ref{levico1} 

\begin{lemma}
\label{su2.1}
Problem \ref{problem1} for $SU(2)$ has no solution if $t\neq 2/k'$ with $k'\in \N$.
\end{lemma}

\begin{proof}
This follows immediately from the proof of lemma \ref{levico1} and from \\
$<\al_1,\al_1>=2$. 
\end{proof}

{}From (\ref{levico2}), (\ref{casca 447}) and (\ref{casca448}) we see that the possible automorphy factors of 
the functions $C^\tau_{2/k'}(\psi) \sigma$ for distributions satisfying (\ref{502}) are
\begin{equation}
\la{maio2}
e(\al + \tau \al', v) = e^{-\pi i k'\al'(v) - \pi i \tau \frac{k'}{2}<\al',\al'>-2\pi i <d,\al'>}, 
\end{equation}
with $d=0$ or $d=(1/2) \lambda_1=(1/4) \al_1$.

{}Functions $\theta^{(1/2)}$ with automorphy factors corresponding to $d=(1/2) \lambda_1$ in (\ref{maio2}) can be obtained 
from those with $d=0$ through a translation of $L_1^{k'}$ by $N_{(\tau \lambda_1/2)}$, 
\begin{equation}
\la{maio3}
\theta(v) \mapsto \theta^{(1/2)}(v) = \theta(v+\lambda_1/k').
\end{equation}
A basis for $H^0(M,N_{(\tau \lambda_1/2)} \otimes L_1^{k'})$ is then given by 
\begin{equation}
\la{maio4}
\theta^{(1/2)}_{m,k'}(z) = \sum_{p\in \Z} (-1)^p e^{\pi i \frac{\tau}{k'}(m+k'p)^2 + 2\pi i (m+k'p)z}, 
\end{equation}
for $m=0,...,k'-1$. The action of the Weyl group $W\cong \Z_2$ is $\theta(v) \mapsto \theta (-v)$ on the elements of the basis 
in (\ref{maio1}) and (\ref{maio4}) is 
\ba
\la{maio5}
\nonumber
\theta_{m,k'}(-z)= \theta_{k'-m,k'} (z) \\
\theta_{m,k'}^{(1/2)}(-z)= \theta_{k'-m,k'}^{(1/2)} (z).
\ea

Let ${\cal H}^+_{k',\tau}$ and  ${\cal H}^-_{k',\tau}$ be the spaces of, respectively, Weyl invariant 
and anti-invariant sections of $L_1^{k'}$. We have the following
\begin{lemma}
\la{maio10}
\begin{itemize}

\item[{1.}] 
If $k'=2k$, with $k\in \N$, then ${\cal H}^+_{2k,\tau}$ has dimension $k+1$ and is generated by 
$\{\theta_{0,2k}, \theta_{k,2k}, \theta_{j,2k}+\theta_{2k-j,2k}, j=1,...,k-1\}$. 
The space ${\cal H}^-_{2k,\tau}$ has dimension 
$k-1$ and a basis is formed by $\{\theta_{j,2k}-\theta_{2k-j,2k}, j=1,...,k-1\}$.
\item[2.] If $k'=2k+1$, with $k\in \N$, then ${\cal H}^+_{2k+1,\tau}$ has dimension $k+1$ and is generated by 
 $\{\theta_{0,2k+1}, \theta_{j,2k+1}+\theta_{2k+1-j,2k+1}, j=1,...,k\}$. 

The space ${\cal H}^-_{2k+1,\tau}$ has dimension 
$k$ and a basis is formed by  $\{\theta_{j,2k+1}-\theta_{2k+1-j,2k+1}, j=1,...,k\}$.
\end{itemize}
\end{lemma}
\begin{proof}
The result follows immediately from (\ref{maio5}). 
\end{proof}

The corresponding similar result holds for the decomposition of  \break 
$H^0(M,N_{(\tau \lambda_1/2)} \otimes L_1^{k'})$ 
under the action of the Weyl group. The analog of theorem \ref{tl} becomes,
\begin{corollary}
\la{maio6} 
\mbox{ }
\begin{itemize}
\item[a)] The space ${\cal H}_{4,\tau}^-$ of Weyl anti-invariant sections of  $L_1^4$ is  
 one-dimensional and ${\cal H}_{4,\tau}^- = <\theta_{4}^->$, where $\theta^-_{4}= \theta_{1,4}-\theta_{3,4}$. 
\item[b)] The map 
\begin{eqnarray}
\label{maio9}
\nonumber
\H_{2k,\tau}^+ & \rightarrow & \H_{2k+4,\tau}^- \\
\theta^+ & \mapsto & \theta^- = \theta_{4}^- \ \theta^+ ,
\end{eqnarray}
is an isomorphism between the space of Weyl invariant sections of $L_1^{2k}$ and the space 
of Weyl anti-invariant sections of $L_1^{2k+4}$. 
\end{itemize}
\end{corollary}
\begin{proof}
This is an immediate consequence of the previous lemma. 
\end{proof}

Again, the corresponding result holds if we replace $L_1^{2k}$ by $N_{(\tau \lambda_1/2)} \otimes L_1^{2k}$. 
We note that the space ${\cal H}_{3,\tau}^-$ is also one-dimensional
and also provides an isomorphism from $\H_{2k,\tau}^+$ to $\H_{2k+3,\tau}^-$, which are both of dimension $k+1$.

Consider the distributions on $SU(2)$ defined by

\begin{eqnarray}
\label{fim2}
\nonumber
\theta^0_{m,k'} (x) & = & \sum_{p\in \Z} e^{2\pi i (m+k'p)x}, \\
\theta^{(1/2) 0}_{m,k'} (x) & = & \sum_{p\in \Z} (-1)^p e^{2\pi i (m+k'p)x},
\end{eqnarray}
for $m=0,...,k'-1$.
We then have,

\begin{theorem}
\la{maio7}
Problem \ref{problem1} has solution if and only if $t= 2/k'$, where $k'$ is of the form 
$k'=2k+3$ or $k'=2k+4$ with $k \in \N_0$. For every such $k'$, 
the set ${\cal F}_k$ of distributions satisfying 
(\ref{502}) is the union of two subspaces ${\cal F}_k^{(0)}$ 
and ${\cal F}_k^{(1/2)}$ intersecting only in $0$ and 
with dimensions equal to $k+1$.
A basis for ${\cal F}_k^{(0)}$ is formed by 
\begin{equation}
\label{fim3}
\psi_{m,k} (v) = \frac{1}{\sigma} \sum_{w\in W} \epsilon(w) \theta^0_{m+1,k'} (w(v)),
\end{equation}
with $m=0,...,k$.
A basis for ${\cal F}_k^{(1/2)}$ is formed by distributions with identical expressions, only with 
$\theta^{(1/2)0}_{m+1,k'}$ replacing $\theta^0_{m+1,k'}$. 
\end{theorem}

\begin{proof}
{}Follows immediately from the proof of theorem \ref{levico6}, with obvious changes 
coming from lemmas \ref{casca449}, \ref{su2.1} and \ref{maio10}.
\end{proof}

When $k'=2k+4$ the image of ${\cal F}_k^{(0)}$ under the CST $C^\tau_{t=2/(2k+4)}$, for 
$\tau \in {\cal T}_1$, gives a vector bundle $\widetilde {\cal H}_{k} \rightarrow {\cal T}_1$, 
which is isomorphic to the bundle of conformal blocks. The analogs of theorems \ref{levico8} and 
\ref{levico9} 
now follow from corollary \ref{maio6} by adapting the proofs of theorems \ref{levico8} and 
\ref{levico9}. 

The image under the CST of the space of distributions ${\cal F}_k^{(1/2)}$, does not lead to pull-backs 
of holomorphic sections of line bundles on ${\cal M}_2 (X_\tau)$. Similarly, 
when $k'=2k+3$, the image by the CST of the space ${\cal F}_k$  leads to holomorphic functions on
$SL(n,\C)$ which are $W \triangleright (\ck \Lambda_R \oplus \tau \ck \lambda_R)$ invariant but which, nevertheless,
descend to ${\cal M}_2(X_\tau)$ only as sections of orbifold line bundles. 

Therefore, while in the case of ${\cal M}_n(X_\tau)$ for $n\geq 3$ problem \ref{problem1}
led automatically to non-abelian theta functions, for $SU(2)$ this happens only for a subset 
of its solutions. 


\section{General Comments and Conclusions}
\label{conclusion}

The vector bundle of conformal blocks $\H_k$ has been the subject of much interest.  
In \cite{AdPW}, Chern-Simons topological quantum field theory was analysed from the point of view of 
geometric quantization on infinite dimensional affine spaces of connections.
In that context, $\H_k$ is endowed with a hermitian structure (see eq. 5.49 in \cite{AdPW}) 
which coincides with the one described in section \ref{s4}. 
In that case, the shift in the level $k \rightarrow k+n$ 
which gives rise to the correct unitary structure on $\H_k$, comes from (not fully rigorous)  
Feynman path integral calculations, leading to  
regularized determinants of Laplace operators on adjoint bundles. On the other hand, 
in the present work the shift follows naturally 
from the generalized CST and the Schottky map. 
In fact, the consideration of class functions or distributions on $K$, 
and the factor of $\sigma$ in the Weyl integration formula, lead to 
Weyl anti-invariance and ultimately to the hermitean structure (\ref{831}) or 
equivalently (\ref{506}). This is the same hermitean structure obtained in 
\cite{AdPW} and, as noted already, is different from the, $a$ $priori$, most natural one, 
for which the frame $\theta^+_{\gamma,k}$ in (\ref{414}) would be orthonormal.
The present work then provides a finite dimensional framework for explaining the shift in the level.

The CST (extended to $(C^\infty(SU(n)))^\prime$ and restricted to ${\cal F}_k$) induces a 
parallel transport on the bundle $\H_k$ which coincides with the one associated to the 
Knizhnik-Zamolodchikov-Bernard 
\cite{KZ,Ber} connection introduced in Wess-Zumino-Witten conformal field theories \cite{Wi1,Wi2}
or also the connection considered in the context of Chern-Simons theories \cite{AdPW,FG,G}. 

We believe that our results also contribute to clarifying the relations between the heat 
equation on the compact group $K$ and 
the representations of the loop group $LK$, a point raised in \cite{PS} (see page 286). 
The genus 1 non-abelian theta functions 
appear in the Weyl-Kac character formula for the affine algebra $\widehat{Lie(LK)}$, 
and as we have seen these theta functions, 
including the unitary structure, can be naturally studied with CST techniques.    

Similar techniques  are also useful for  the (much harder) study of non-abelian theta functions for curves 
of genus greater than 1 \cite{FMNT}. 
We note that in the present work and also in the case of classical (abelian) theta functions \cite{FMN}, the choice of 
distributions to which the CST should be applied to produce non-abelian theta 
functions is dictated by naturality conditions 
related to the unitarity of the extended CST. Alternatively, in both cases 
the same distributions can be determined by   
Bohr-Sommerfeld conditions in geometric quantization \cite{Ty}. 
In fact, the distributions in (\ref{510}) are combinations of Dirac delta distributions 
supported on Bohr-Sommerfeld points. This should be related to the work of \cite{We}. 
We expect that this continues to be true in the higher genus non-abelian case as well.


\section{Appendix}
\la{app}

In this appendix, for convenience of the reader, we obtain explicit
expressions for the canonical bases of $\check{\Lambda}_{R}\oplus \tau 
\check{\Lambda}_{R}$ with respect to the form $E_{1}$ in (\ref{411}), 
for the case
of the Lie algebra $sl(n,\Bbb{C}).$ We also write down the corresponding
period matrices.

We start with a lemma valid for a simple Lie algebra $\mathfrak g$,
of rank $l$. To simplify the notation, instead of the coroot lattice, 
we will work with the root lattice $\Lambda _{R}\subset\mathfrak{h}$.
Let $E$ be the antisymmetric $\Bbb{Z}$-bilinear form on 
$\Lambda _{R}\oplus \tau \Lambda _{R}$ given by 
\[
E\left( \alpha _{i},\alpha _{j}\right) =E\left( \tau \alpha _{i},\tau \alpha
_{j}\right) =0,\qquad E\left( \alpha _{i},\tau \alpha _{j}\right) =-E\left(
\tau \alpha _{j},\alpha _{i}\right) =C_{ij}, 
\]
where $C=\left[ C_{ij}\right] $ is the Cartan matrix of $\mathfrak g$. Recall that a
canonical basis  $(\beta _{1},...,\beta _{l},\tilde{\beta}_{1},....,\tilde{%
\beta}_{l})$ of $\Lambda _{R}\oplus \tau \Lambda _{R}$ with respect to $E$
satisfies by definition 
\[
E(\beta _{i},\beta _{j})=E(\tilde{\beta}_{i},\tilde{\beta}_{j})=0,\qquad
E(\beta _{i},\tilde{\beta}_{j})=-E(\tilde{\beta}_{i},\beta _{j})=\delta
_{i}\delta _{ij}, 
\]
where $\delta _{1}|\cdots|\delta _{l}$ are integers, 
with $\delta_1\cdots\delta_l=\det C$.

To find such a canonical basis, let us write the (ordered) basis of simple roots in $%
\Lambda _{R}$ as a row vector ${\bf\alpha }=\left( \alpha
_{1},...,\alpha _{l}\right) $. In particular, any other basis ${\bf\beta 
}=\left( \beta _{1},...,\beta _{l}\right) $ can be written as ${\bf\beta 
}={\bf\alpha }A$ for some integer matrix $A$ with unit determinant. In
this notation, it is easy to see that the pair $({\bf\beta },{\bf%
\tilde{\beta}})$ forms a canonical basis of $\Lambda _{R}\oplus \tau \Lambda
_{R}$ if and only if there are matrices $A,\tilde{A}\in SL(l,\Bbb{Z})$ such
that 
\begin{equation}
A^{t}C\tilde{A}=\Delta ,  \tag{*}
\end{equation}
where $\Delta =diag(\delta_1,...,\delta_n),$ 
and $A^{t}$ denotes the transpose of $A$.
In fact, the bases are related by ${\bf\beta }={\bf\alpha }A$ and $%
{\bf\tilde{\beta}}=\tau {\bf\alpha }\tilde{A}.$

The solutions of (*) can be naturally identified with bases of $\Lambda
_{W}, $ the weight lattice. More precisely, let us say that a basis ${\bf%
\beta }$ of $\Lambda _{R}$ is {\em completable} 
if there exists a basis $({\bf%
\beta },{\bf\tilde{\beta}})$ of $\Lambda _{R}\oplus \tau \Lambda _{R}$,
which is canonical with respect to $E,$ or equivalentely, if there exists a
solution of (*) with ${\bf\beta }={\bf\alpha }A.$ Then, we have

\begin{lemma}
${\bf\beta }$ is completable if and only if ${\bf\beta }\Delta ^{-1}$
is a basis of $\Lambda _{W}.$
\end{lemma}

\begin{proof}
Let ${\bf\beta }={\bf\alpha }A$ and let $A^{t}C\tilde{A}%
=\Delta $ for some matrices $A,\tilde{A}\in SL(l,\Bbb{Z}).$ Since the
relation between roots and weights is given by ${\bf\alpha }={\bf%
\lambda }C,$ we have 
\[
{\bf\beta }\Delta ^{-1}={\bf\alpha }A\Delta ^{-1}={\bf\lambda }%
CA\Delta ^{-1}, 
\]
which means that ${\bf\beta }\Delta ^{-1}$ is a basis of $\Lambda _{W}$,
because $CA\Delta ^{-1}=( \tilde{A}^{t}) ^{-1}$ is a unimodular
matrix.
\end{proof}

Now let us consider the case of $sl(n,\Bbb{C})$, 
where the simple roots $\left\{ \alpha
_{1},...,\alpha _{l}\right\} $ form a basis of $\Lambda _{R}$, $n=l+1$. 
From the above lemma, we see that ${\bf\alpha }$ is not a completable basis
of $\Lambda _{R},$ except when $l=1$ (in which case $\left( {\bf\alpha }%
,\tau {\bf\alpha }\right) $ is clearly a canonical basis of
$\Lambda _{R}\oplus \tau \Lambda _{R}$). 
We assume henceforth that $l\geq 2,$ and let us
consider a slightly different basis 
\[
{\bf\beta =(}\beta _{1},...,\beta _{l}):=(\alpha _{1},...,\alpha
_{l-1},n\lambda _{1}), 
\]
where $\lambda_1$ is the fundamental weight dual to $\alpha_1$.
Since $n\lambda _{1}=l\alpha _{1}+\left( l-1\right) \alpha _{2}+...+\alpha
_{l},$ then ${\bf\beta }={\bf\alpha }A$ for a unimodular matrix $A,$
which means that ${\bf\beta }$ is another basis of $\Lambda _{R}.$ 
To prove that this ${\bf\beta }$ 
is indeed completable, consider
$\delta _{1}=...=\delta _{l-1}=1$ and $\delta _{l}=n$,
and write ${\bf\beta }\Delta
^{-1}=(\alpha _{1},...,\alpha _{l-1},\lambda _{1})$; 
this implies that the matrix
relating ${\bf\lambda }$ to ${\bf\beta }\Delta ^{-1}$ is an integer
matrix; since this matrix is $CA\Delta ^{-1}$ (from the lemma), it has
determinant one, which proves that ${\bf\beta }\Delta ^{-1}$ is a
basis of $\Lambda _{W}.$ (Note that $\det C=\det \Delta =n).$

Simple calculations now give the following explicit expressions.
Recall that the period matrix $\Omega$ is defined by 
${\bf\tilde{\beta}}={\bf\beta }\Omega \Delta $.

\begin{proposition}
For $l\geq 2,$ a canonical basis of $\Lambda _{R}\oplus \tau \Lambda _{R}$
with respect to $E$ is $(\beta _{1},...,\beta _{l},\tilde{\beta}_{1},....,%
\tilde{\beta}_{l}),$ where 
\[
\left\{ 
\begin{array}{ll}
\beta _{i}=\alpha _{i}, & i=1,...,l-1 \\ 
\beta _{l}=n\lambda _{1}=l\alpha _{1}+\left( l-1\right) \alpha
_{2}+...+\alpha _{l} &  \\ 
\tilde{\beta}_{i}=\tau \lambda _{i}-\left( n-i\right) \tau \lambda _{l}, & 
i=1,...,l-1 \\ 
\tilde{\beta}_{l}=n\tau \lambda _{l} & 
\end{array}
\right. 
\]
The period matrix is then given by 
\[
\left\{ 
\begin{array}{ll}
\Omega _{i,i+k}=\Omega _{i+k,i}=\tau (n-i)(l-i-k), & i+k\leq l-1 \\ 
\Omega _{i,l}=\tau (i-l), & i\leq l-1 \\ 
\Omega _{l,l}=\tau \frac{l}{n}. & 
\end{array}
\right. 
\]
\end{proposition}

\begin{proof} 
We have verified the formula for ${\bf\beta ,}$ and the one for $%
{\bf\tilde{\beta}}=(\tilde{\beta}_{1},....,\tilde{\beta}_{l})$ is a
direct calculation from ${\bf\tilde{\beta}}=\tau {\bf\alpha }\tilde{A%
}=\tau {\bf\lambda }\left( A^{t}\right) ^{-1}\Delta $ (using (*)). 
Since we also have ${\bf\tilde{\beta}}=\tau {\bf\alpha }\tilde{A}%
=\tau {\bf\beta }A^{-1}\tilde{A},$ we get 
\[
\Omega =\tau A^{-1}\tilde{A}\Delta ^{-1}=\tau A^{-1}C^{-1}(A^{t}) ^{-1},
\]
which readily gives the above expression for the period matrix (note that
all its entries have to lie in $\frac{\tau}{n}\Bbb{Z}).$
\end{proof}

Let us now examine the other possible canonical bases, i.e, all the
solutions to (*). If the pairs $(A_{1},\tilde{A}_{1})$ and $(A_{2},\tilde{A}%
_{2})$ are both solutions of (*), with $A_{2}=A_{1}B,$ and $\tilde{A}_{2}=%
\tilde{A}_{1}\tilde{B},$ for some unimodular matrices $B$ and $\tilde{B}$,
then necessarily 
\begin{equation}
B^{t}\Delta \tilde{B}=\Delta .  \tag{**}
\end{equation}

It is not difficult to see that the set of matrices $B$ such that there
exists $\tilde{B}$ satisfying (**) form a group, and that this is the
subgroup $\Gamma _{n}\subset SL(l,\Bbb{Z})$ of modular matrices of the form 
\[
B=\left( 
\begin{array}{cc}
a & b \\ 
c & d
\end{array}
\right) 
\]
where $a$ is an $\left( l-1\right) \times \left( l-1\right) $ matrix, $d$ is
an integer, $b$ and $c^{t}$ are $\left( l-1\right) $-vectors and 
$b\in n \Bbb{Z}^{l-1}$ (all entries in $b$ are multiples of $n$). Finally, if 
$\Omega _{i}$ are period matrices given by $\Omega _{i}=\tau
A_{i}^{-1}C^{-1}\left( A_{i}^{t}\right) ^{-1},$ $i=1,2,$ as in the
proposition above, and $A_{2}=A_{1}B^{-1},$ for some $B\in \Gamma _{n},$
then we have 
\[
\Omega _{2}=B\Omega _{1}B^{t}.
\]

This means that two period matrices which are related in this way for some 
$B\in \Gamma _{n}$ should be regarded as equivalent.

\section*{Acknowledgements:}
We wish to  thank Professor M.S. Narasimhan
for discussions and suggestions at an early stage
of this work. It is a pleasure  to thank Professor Andrei Tyurin for 
many fruitful discussions and encouragement. 
The authors were partially supported by 
the Center for Applied Mathematics, IST, Lisbon (CF and JM),   by 
the Center for Analysis, Geometry and Dynamical Systems, IST, Lisbon (JPN),  
and also by the FCT (Portugal) via the program POCTI and by the projects POCTI/33943/MAT/2000 
and CERN/FIS/43717/2001.





\end{document}